\input amstex
\documentstyle{amsppt}
\pageno=1
\magnification1200
\catcode`\@=11
\def\logo@{}
\catcode`\@=\active
\NoBlackBoxes
\vsize=23.5truecm
\hsize=16.5truecm

\def\d{d\!@!@!@!@!@!{}^{@!@!\text{\rm--}}\!}

\def\Zfrac{\tsize\frac1{\raise 1pt\hbox{$\scriptstyle z$}}}

\def\crp{\overline{\Bbb R}_+}
\def\crm{\overline{\Bbb R}_-}
\def\crpm{\overline{\Bbb R}_\pm}

\def\rnp{{\Bbb R}^n_+}

\def\ang#1{\langle {#1} \rangle}

\def\ttilde{\overset{\,\approx}\to}

\def\rp{ \Bbb R_+}

\define\tr{\operatorname{tr}}

\define\Tr{\operatorname{Tr}}

\define\srpp{\Cal S_{++}}
\define\stimes{\!\times\!}

\document

\topmatter
\title
{Spectral boundary conditions for generalizations of Laplace and
Dirac operators}  
\endtitle
\author{Gerd Grubb}
\endauthor
\affil
{Copenhagen Univ\. Math\. Dept\.,
Universitetsparken 5, DK-2100 Copenhagen, Denmark.
E-mail {\tt grubb\@math.ku.dk}}\endaffil
\rightheadtext{Spectral boundary conditions}

\abstract
Spectral boundary conditions for Laplace-type operators on a compact
manifold $X$ with boundary are partly Dirichlet, partly (oblique)
Neumann conditions, where the partitioning is provided by a
pseudodifferential projection; they have an interest in string and
brane theory. 
Relying on pseudodifferential methods, we give sufficient conditions
for the existence of the associated 
resolvent and heat operator, and show asymptotic expansions of their
traces in powers and power-log terms, allowing a smearing function
$\varphi $. The leading log-coefficient is identified as a
non-commutative residue, which vanishes when $\varphi =1$.

The study has new consequences for well-posed (spectral) boundary problems for
first-order, Dirac-like elliptic operators (generalizing the
Atiyah-Patodi-Singer problem). 
It is found e.g\. that the zeta function is always regular at zero. In the
selfadjoint case, there is a stability of the zeta function value and
the eta function regularity at zero, under perturbations of the
boundary projection of order $-\operatorname{dim}X$.
\endabstract

\endtopmatter

\subhead 
\footnote{\eightrm To appear in Comm\. Math\. Phys.}Introduction
\endsubhead

Spectral boundary conditions (involving
a pseudodifferential projection $\Pi $ on the boundary)
were first employed by Atiyah, Patodi and Singer \cite{APS75} in their
seminal work on Dirac operators $D$ on manifolds $X$ with boundary,
introducing the eta-invariant of the tangential part of $D$ as an
extremely interesting new geometric object. 
For such a Dirac realization $D_\Pi $, the operator ${D_\Pi}^* D_\Pi $
(the square ${D_\Pi }^2$ in the selfadjoint case) is  a
Laplace operator with a boundary condition of the form
$$
\Pi \gamma _0u=0,\quad (I-\Pi )(\gamma _1u+B\gamma _0u)=0;
\tag0.1
$$
here $\gamma _ju=(\frac \partial {\partial n})^ju|_{\partial X}$, and
$B$ denotes a first-order operator on $\partial X$ (in the case
derived from $D_\Pi $ it is the tangential part of $D$). 
In this work we present a new analysis of Laplace-type
operators with boundary conditions (0.1), under quite general choices
of $\Pi $ and $B$, showing existence of heat trace expansions and
meromorphic zeta functions, and analyzing the leading logarithmic and
nonlocal terms. The methods developed here moreover lead to new
results on the corresponding questions for Dirac operator problems.

In a physics context, spectral boundary conditions are used in 
studies of axial anomalies (see e.g.\ Hortacsu, Rothe and Schroer
\cite{HRS80}, Ninomiya and Tan \cite{NT85}, Niemi
and Semenoff \cite{NS86}, Forgacs, O'Raifertaigh and Wipf
\cite{FOW87}, see also 
Eguchi, Gilkey and Hanson \cite{EGH80}), 
in quantum cosmology (see e.g.\ D'Eath and Esposito \cite{EE91}), and in the theory
of the Aharonov-Bohm effect (see e.g.\ Beneventano, De Francia and
Santangelo \cite{BFS99}).
For a long time, applications of spectral boundary conditions
were limited to fields of half-integer spin whose dynamics
are governed by first-order operators of Dirac type. 
It is natural (and even required by supersymmetry arguments)
to extend this scheme to integer spin fields, and, therefore,
to operators of Laplace type. Indeed, a first step in this
direction has been taken in Vassilevich \cite{V01}, \cite{V02},
where spectral boundary conditions are formulated for bosonic
strings to describe some collective states of open strings and
Dirichlet branes.

Local cases of (0.1), where $\Pi $ is a projection morphism and $B$ is
a differential operator, have been
treated earlier (see e.g\. Avramidi and Esposito \cite{AE99} and its
references); such cases have a
mathematical foundation in  Greiner \cite{Gre71}, Grubb \cite{G74},
Gilkey and Smith \cite{GiS83}, whereas the general global case is a new subject.

Establishing heat trace asymptotics is very important in quantum field theory
since they are related to ultra violet divergences and quantum anomalies
(see e.g.\ the surveys of Avramidi \cite{A02} and Fursaev \cite{F02}). A vanishing
of the leading logarithmic term means that the standard definition of
a functional 
determinant through the zeta function derivative at 0 is applicable,
and the road is open to renormalization  
of the effective action.

\smallskip
\flushpar{\it Overview of the contents:} 
The problems for Laplace-type operators $P$ are considered in
Sections 1--4. 
In Section 2 we give sufficient conditions for the
existence of the resolvent $(P_T-\lambda )^{-1}$ in an
angular region, with 
an explicit formula (Th\. 2.10), and we use this to show asymptotic
expansions of traces $\Tr(F (P_T-\lambda )^{-m})$ (for sufficiently large
$m$) in powers and log-powers of $\lambda $ (Th\. 2.13); here $F$ is an arbitrary differential operator. In Section 3
we show corresponding expansions of the heat trace $\Tr(Fe^{-tP_T})$
in powers and log-powers of $t$ (Cor\. 3.1), and 
establish meromorphic
extensions of zeta functions $\zeta (F,P_T,s)=\Tr(FP_T^{-s})$ 
with simple and double real poles in $s$. (These constructions do not
need a differential operator square root of $P$ as discussed in \cite{V01}.)

The leading logarithmic term and
nonlocal term are analyzed in Section 4 when $F$ is a morphism
or a first-order operator. This is done by determining which part of the
resolvent actually contributes to these values (Th\. 4.1).
It follows that the
log-coefficient identifies with a certain
non-commutative residue, vanishing e.g.\ when $F=I$ (Th\. 4.2, 4.5).
Then the zeta function is regular
at 0, and the value at zero can be set in relation to an eta-invariant
associated with $\Pi $ (Th\. 4.9, Def\. 4.10). Special results
are also obtained when $F$ is a first-order operator and certain symmetry
properties hold.

In Section 5 we draw some conclusions
for operators 
${D_\Pi }^*D_\Pi $ defined from a Dirac-type operator $D$ (of the
form $\sigma (\partial _{x_n}+A+\text{perturb.})$ near $\partial X$
with $A$ selfadjoint elliptic on $\partial X$) together
with a pseudodifferential orthogonal projection $\Pi $ on $\partial X$
defining the boundary condition $\Pi \gamma _0u=0$.
Atiyah, Patodi and Singer \cite{APS75} considered such problems with $\Pi
$ equal to the nonnegative eigenprojection $\Pi _\ge(A)$, and
increasingly general choices have been 
treated through the years: Douglas and Wojciechowski \cite{DW91},
M\"uller \cite{M94}, 
Dai and Freed \cite{DF94},
Grubb and Seeley \cite{GS95}, \cite{GS96} 
gave results on perturbations \linebreak$\Pi _\ge(A)+\Cal S$ by finite rank
operators $\Cal S$; \cite{W99} allowed smoothing operators $\Cal S$.  
Br\"uning and Lesch \cite{BL99} introduced a special class of other
projections $\Pi (\theta )$, and \cite{G99}, \cite{G01$'$} included all
projections satisfying the well-posedness condition of Seeley \cite{S69$''$}. 

The methods of the present study lead to new results both for zeta and
eta functions. The zeta function $\zeta (D_\Pi ^*D_\Pi ,s)=\Tr ((D_\Pi
^*D_\Pi )^{-s}) $ is shown to be regular at
$s=0$ for {\it any
well-posed} $\Pi $ such that
the principal parts of $\Pi $
and $A^2$ commute (Cor\. 5.3). This was known previously for
perturbations of $\Pi _\ge(A)$ and $\Pi (\theta )$ of order $-\dim X$
(cf\. \cite{G01$'$}); the new result
includes in particular the perturbations of $\Pi _\ge(A)$ and $\Pi
(\theta )$ of order
$-1$. The value $\zeta (D_\Pi ^*D_\Pi ,0)$ is determined from
$\Pi $ and 
$\operatorname{ker}D_{\Pi }$ modulo local terms (Cor\. 5.4).

Restricting the attention to cases with selfadjointness at $\partial X$,
assuming $$
\sigma A=-A\sigma ,\quad \sigma ^2=-I,\quad \Pi =-\sigma \Pi
^\perp\sigma ,\tag0.2
$$ we get new results for the eta function $\eta
(D_\Pi ,s)=\Tr (D(D_{\Pi }^*D_\Pi )^{-\frac{s+1}2})$: It has at most
a simple pole at 
$s=0$ for general well-posed  $\Pi $ (Cor\. 5.8); this includes
perturbations of  
$\Pi _\ge(A)$ and $\Pi
(\theta )$ of order
$-1$, where it was previously known for order $-\dim X$ in the
selfadjoint product 
case. Moreover, the residue at $s=0$ is locally determined, and stable
under perturbations of 
$\Pi $ of order $\le -\dim X$ (Th\. 5.9). In particular, in the
selfadjoint product case, the vanishing of the simple pole, shown for special
cases in \cite{DW91}, \cite{M94}, \cite{W99}, \cite{BL99}, is stable under perturbations
of order $-\operatorname{dim}X$. (For the last result, see also Lei \cite{L02}.)
There are similar stability results for the value of zeta at 0 (Th\. 5.7).

\smallskip
The author is grateful to D\. Vassilevich, R\. Mazzeo, P\.
Gilkey, R\. Melrose, E\. Schrohe and B\. Booss-Bavnbek for useful
conversations.

\subhead 1. Boundary conditions with projections
\endsubhead

Consider  a second-order strongly elliptic
differential operator $P$ acting on the sections of an
$N$-dimensional $C^\infty $ vector bundle $E$ over a compact $C^\infty $
$n$-dimensional manifold $X$ with boundary $X'=\partial X$.

$X$ is provided with a volume element and $E$ with a hermitian
metric defining a Hilbert space structure on the sections, $L_2(E)$.
We denote $E|_{X'}=E'$. A neighborhood of
$X'$ in $X$ has the form $X_c=X'\times[0,c[$ (the points denoted
$x=(x',x_n)$), and there $E$
is isomorphic to the pull-back of $E'$. On $X_c$, there is a smooth
volume element $v(x)dx'dx_n$, and
$v(x',0)dx'$ is the volume element on $X'$ (defining $L_2(E')$).

We assume that $P$ is principally selfadjoint, i.e., $P-P^*$ is of
order $\le 1$. Moreover, to have simple ingredients to work with,
we assume that
$P$ is of the following form near $X'$:

\proclaim{Assumption 1.1} On $X_c$, $P$ has the form$$
P=-\partial _{x_n}^2+P'+x_nP_2+P_1,\tag1.1
$$ 
where $P'$ is an elliptic selfadjoint
nonnegative second-order
differential operator in $E'$ (independent of $x_n$) and the $P_j$ are
differential operators of order $j$ in $E|_{X_c}$.  
\endproclaim 

Let $\Pi _1$  be a classical pseudodifferential operator ($\psi $do)
in $E'$ of 
order 0 with $\Pi _1^2=\Pi _1$, i.e., a projection operator, and
denote the complementing projection by $\Pi _2$:$$
\Pi _2=I-\Pi _1.\tag1.2
$$ 
Let
$B$ be a first-order differential or pseudodifferential operator in
$E'$. 
Then we consider the boundary condition for $P$:$$
\Pi _1\gamma _0u=0,\quad \Pi _2(\gamma _1u+B\gamma _0u)=0,\tag1.3
$$
where the notation $\gamma _ju=(\partial _{x_n}^ju)|_{X'}$ is used.
In short, $Tu=0$, where 
$$T=\{\Pi _1\gamma _0,\Pi _2(\gamma _1+B\gamma
_0)\}.\tag1.4
$$
We shall study the resolvent and the heat operator defined from
$P $ under this boundary condition, when suitable
parameter-ellipticity conditions are satisfied.
More precisely, with $H^s(E)$ denoting the Sobolev
space of order $s$ (with norm $\|u\|_s$), we define the realization
$P_T$ in $L_2(E)$ determined by the boundary condition
(1.3) as the operator acting like $P$ and with domain 
$$
D(P_T)=\{u\in H^2(E)\mid Tu=0\};
$$ then we
want to construct the resolvent $(P_T-\lambda )^{-1}$ and the heat
operator $e^{-tP_T}$ and analyze their trace properties. In
particular, we want to show a heat trace expansion
$$
\Tr  e^{-tP_T}
\sim
\sum_{-n\le k<0}  a_{ k}t ^{\frac{ k}2}+ 
\sum_{k\ge 0}\bigl({  -a'_{ k}}\log t +{ a''_{
k}}\bigr)t^{\frac{ k}2}.\tag1.5
$$
Such expansions are known to hold for normal differential boundary
conditions (Seeley \cite{S69}, \cite{S69$'$}, Greiner \cite{Gre71}) without the
logarithmic terms, and for pseudo-normal
$\psi $do boundary conditions \cite{G99}, but the condition (1.3) is in
general not of these types.

\example{Example 1.2} Assumption 1.1 holds if $P=D^*D$ for an
operator type including Dirac operators, namely a first-order elliptic differential operator
$D$ from $E$ to another $N$-dimensional bundle $E_1$ over $X$ such that$$
D=\sigma (\partial _{x_n}+A_1)\;\text{ with }\; A_1=A+x_nA_{11}+A_{10}\;\text{ on }X_c,\tag1.6 
$$ 
as considered
in numerous works, originating in Atiyah, Patodi and Singer \cite{APS75};
here $\sigma $ is a unitary morphism from $E'$ to $E'_1$, $A$ is a selfadjoint first-order elliptic operator in $E'$
independent of $x_n$, and the $A_{1j}$ are $x_n$-dependent
differential operators in $E'$ of order $j$. Then we can take
$P'=A^2$. However, for   
the study in the following one does not need to have a differential
operator ``square root'' 
of $P'$. (And, if $P$ is derived from
(1.6) near
$X'$, the factorization need not extend to all of $X$.) 
See Vassilevich \cite{V01} for a discussion of how one can find such $D$ when
$P$ is a Laplacian.

When $P=D^*D$, one can
for example take as 
$\Pi _1$ the orthogonal 
projection $\Pi _\ge(A)$ onto the 
nonnegative eigenspace of $A$ and let
$B=A_1|_{x_n=0}$, also denoted $A_1(0)$. Then $P_T$ equals $D^*_\ge D_\ge$, where $D_\ge$ is the
realization of $D$ under the boundary condition $\Pi _\ge(A)\gamma
_0u=0$, and (1.5) is known from \cite{GS95} (also with certain finite rank
perturbations of $\Pi _\ge(A)$). 
Boundary conditions (1.3) with such choices of $\Pi _1$ are often
called {\it 
spectral boundary conditions}.

The paper \cite{G99} allows much more general realizations $D_{\Pi _1}$ of
$D$, where $\Pi _\ge(A)$
is replaced by a $\psi $do projection $\Pi _1$ that is ``well-posed'' with
respect to $D$. Well-posedness (introduced by Seeley \cite{S69$''$}) means that the principal symbol $\pi _1^0$ at each
$(x',\xi ')\in S^*(X')$ maps $N_+(x',\xi ')$ bijectively onto the
range of $\pi _1^0(x',\xi ')$ in $\Bbb C^N$, where $N_+(x',\xi ')$ is
the space of boundary values of null-solutions:
$$
N_+(x',\xi ')=\{\, z(0)\in\Bbb C^N\mid d^0(x',0,\xi ',D_{x_n})z(x_n)=0,\, z\in
L_2(\Bbb R_+)^N\,\},
$$
$d^0$ denoting the principal symbol of $D$. As shown in \cite{S69$''$},
the closed-range operator defining the boundary condition may always
be taken to be a projection; it may even be taken orthogonal.
When $\Pi _1$ is orthogonal, then $(D_{\Pi
_1})^*D_{\Pi _1}=P_T$ with $B=A_1(0)$ in 
(1.3)--(1.4), and
(1.5) holds. (Further details in Section 5.) 

In \cite{V01}, Vassilevich points to the need for considering other choices of
$B$ --- and not just $B=0$ but moreover cases unrelated to $A_1$ and having
complex coefficients, where 
$P_T$ is not selfadjoint. He inquires about heat kernel results
for such cases, and that is precisely what we want to develop here.
\endexample

The conditions in (1.3) can be defined without reference to a
factorization as in Example 1.2;
$\Pi _1$ can be a general pseudodifferential projection in $E'$
unrelated to $P'$.

For the case where $\Pi _1$ is a local projection
(a projection 
morphism in $E'$ ranging in a subbundle $F_0$), the question of heat trace expansions is covered
by Greiner \cite{Gre71}, see also \cite{G74}, \cite{G96}.

Since $P$ is strongly elliptic, it satisfies the G\aa{}rding
inequality, so we can assume (possibly after addition of a constant) that
$$
\operatorname{Re}(Pu,u)\ge c_0\|u\|_1^2, \text{ for }u\in C_0^\infty
(X^\circ), \tag1.7
$$
with $c_0>0$.
Then the Dirichlet realization $P_{\operatorname{D}} $ (the realization
of $P$ in $L_2(E)$ with domain $D(P_{\operatorname{D}} )=\{u\in H^2(E)\mid \gamma _0u=0\}$) is
invertible. We shall use the following notation for regions in $\Bbb
C$:
$$
\Gamma _{\theta }=\{\, \mu \in\Bbb C\setminus \{0\}\mid |\arg\mu
|<\theta\,\},\quad \Gamma =\Gamma _{\frac\pi 2},\quad \Gamma _{\theta ,r}=\{\, \mu \in\Gamma _\theta
\mid |\mu |>r\,\}.\tag1.8
$$
Since the principal symbol of $P$ is positive
selfadjoint, the spectrum of $P_{\operatorname{D}} $ is for any $\delta
>0$ contained in 
a set $$
\Sigma _{\delta ,R}=
\Gamma _{\delta }\cup\{|\lambda |\le R(\delta )\} \tag1.9$$ for some
$R=R(\delta )$,
in addition to being contained in $\{\operatorname{Re}\lambda >0\}$.

We denote $(P_{\operatorname{D}}-\lambda
)^{-1}=R_{\operatorname{D}}(\lambda )$, the resolvent of the
Dirichlet problem, i.e\. the solution operator for the
semi-homogeneous problem $$
(P-\lambda )u=f\text{ in }X,\quad \gamma _0u=0 \text{ on }X'; \tag1.10
$$ 
$R_{\operatorname{D}}$ maps $H^s(E)$ continuously into $H^{s+2}(E)$
for $s>-\frac32$, when
$\lambda \notin \Sigma _{\delta , R}$. The other semi-homogeneous
Dirichlet problem 
$$
(P-\lambda )u=0\text{ in }X,\quad \gamma _0u=\varphi  \text{ on }X', \tag1.11
$$ 
is likewise uniquely solvable for $\lambda \notin\Sigma _{\delta ,R}$;
the solution operator will be denoted $K_{\operatorname{D}}(\lambda
)$. This is 
an elementary Poisson operator (in the notation of Boutet de Monvel
\cite{BM71}), mapping $H^s(E')$ continuously into $H^{s+\frac12}(E)$ for
all $s\in\Bbb R$.

We have Green's formula$$
(Pu,v)_{X}-(u,P^*v)_{X}=
(\gamma _1u+\sigma _0\gamma _0u,\gamma _0v)_{X'}-(\gamma
_0u,\gamma _1v)_{X'},\tag1.12
$$
with a certain morphism $\sigma _0$ in $E'$.

\subhead 2. Resolvent constructions \endsubhead

Our way to construct the heat operator for $P_T$ goes via the
resolvent, which is particularly well suited to calculations in the
pseudodifferential framework, since the spectral parameter $\lambda  $
(or,
rather, the square root $\mu =\sqrt{-\lambda }\,$) enters to some
extent like a cotangent variable.
In the following, we shall freely use the notation and results of
Grubb and Seeley \cite{GS95}, Grubb \cite{G01}. To save space, we do not repeat
many details here but 
refer to these 
papers or to the perhaps simpler resum\'e of the needed parts of the
calculus in \cite{G02, Sect\. 2}. Let us just recall the definition of the
symbol space $S^{m,d,s}(\Bbb R^{n-1}\stimes\Bbb R^{n-1},\Lambda )$
(denoted $S^{m,d,s}(\Lambda )$ for short), where $\Lambda $ is a
sector of $\Bbb C\setminus\{0\}$ and $m,d,s\in\Bbb Z$:

A $C^\infty $ function $p(x',\xi ',\mu )$ lies in $S^{m,0,0}(\Bbb R^{n-1}\stimes\Bbb
R^{n-1},\Lambda )$ when, for every closed 
subsector $\Lambda '$,
$$\partial _{|z|}^jp(x',\xi ',1/z)\in S^{m +j}(\Bbb R^{n-1}\stimes\Bbb
R^{n-1})\text{ uniformly for }|z|\le 1,1/z\in \Lambda ';
$$
here $S^{k}(\Bbb R^{n-1}\stimes\Bbb R^{n-1})$ is the usual symbol space of
functions $q(x',\xi ')$ with $$
|\partial _{x'}^\beta \partial _{\xi
'}^\alpha q(x',\xi ')|\le C_{\alpha ,\beta }\ang{\xi '}^{k-|\alpha |}$$
for all $\alpha ,\beta \in\Bbb N^{n-1}$; we write $\ang
x=(1+|x|^2)^{\frac12}$ and  $\Bbb N=\{0,1,2,\dots\}$. Moreover,
$$
S^{m,d,s}(\Bbb R^{n-1}\stimes\Bbb R^{n-1},\Lambda ) =
\mu ^d(|\mu |^2+|\xi '|^2)^{s/2}S^{m,0,0}(\Bbb R^{n-1}\stimes\Bbb R^{n-1},\Lambda ).
$$
In the applications, we often need $p$ and its derivatives to be
holomorphic in $\mu \in \Lambda ^\circ$ for $|\mu |+|\xi '|\ge
\varepsilon >0$;  such
symbols will just be said to be holomorphic (in $\mu $).

The space denoted $S^{m,d}(\Lambda )$ in \cite{GS95} is the space of
holomorphic symbols in $S^{m,d,0}(\Lambda )$. The third upper index $s$
was added in \cite{G01} for convenience; one has in view of \cite{GS95, Lemma
1.13} that$$
\aligned
S^{m,d,s}(\Lambda )&\subset S^{m+s,d,0}(\Lambda )\cap
S^{m,d+s,0}(\Lambda )\text{ for }s\le 0,\\
S^{m,d,s}(\Lambda )&\subset S^{m+s,d,0}(\Lambda )+
S^{m,d+s,0}(\Lambda )\text{ for }s\ge 0,
\endaligned\tag2.1
$$
and the $s$-index saves us from keeping track of a lot of sums and
intersections.

The $\psi $do with symbol $p(x',\xi ',\mu )$ is defined by the usual
formula:
$$
\operatorname{OP}'(p)\:v(x')\mapsto \int_{\Bbb
R^{2(n-1)}}e^{i(x'-y')\cdot
\xi '} p(x',\xi ',\mu )v(y')\,dy'\d\xi ',
$$
where
$\d \xi '$ stands for $(2\pi )^{1-n}d\xi '$; the analogous
definition on $\Bbb R^n$ is indicated by OP. $\Psi $do's in bundles
over manifolds are defined by use of local trivializations.

The symbols, we consider, moreover have expansions in
homogeneous terms, 
$p\sim$ $\sum_{j\in\Bbb N}p_{m-j}$ with  $p_{m-j}\in
S^{m-j,d,s}(\Lambda)$, homogeneous in $(\xi ',\mu )$ of degree
$m-j+d+s$. In the general, so-called weakly polyhomogeneous case, the
homogeneity  takes place for $|\xi '|\ge 1$, 
but if it extends to $|\xi '|+|\mu |\ge 1$ (in such a way that the symbol
behaves as a standard classical symbol in the non-parametrized calculus with
an extra cotangent variable $|\mu |$ entering 
on a par with $\xi '$ in the estimates), the
symbol is
called strongly polyhomogeneous.

The composition rules for these spaces are straightforward: When
$p_i\in S^{m_i,d_i,s_i}(\Lambda )$ for $i=1,2$, then 
$p_1p_2$ and $p_1\circ p_2\in S^{m_1+m_2,d_1+d_2,s_1+s_2}(\Lambda )$;
the latter is the symbol of the composed operator
$\operatorname{OP}'(p_1)\operatorname{OP}'(p_2)$, satisfying 
$$
(p_1\circ p_2)(x',\xi ',\mu )\sim\sum_{\alpha \in\Bbb 
N^{n-1}}\tfrac {(-i)^{|\alpha |}}{\alpha
!}\partial _{\xi '}^\alpha p_1 (x',\xi ',\mu )\partial _{x'}^\alpha p_2(x',\xi
',\mu ).\tag 2.2
$$

Before discussing the construction of the resolvent, we shall introduce
some auxiliary pseudodifferential operators on $X'$. When $\lambda \in \Bbb
C\setminus\rp$, we write $\mu =(-\lambda)^\frac12
$, defined such that $\mu \in \Gamma $ (cf.\ (1.8)).

\proclaim{Definition 2.1} The operator $\frak A(\lambda )$ is defined
for $\lambda \in \Bbb C\setminus\rp$ by$$
\frak A(\lambda )=(P'-\lambda )^\frac12;
$$
it is a $\psi $do in $E'$ of order 1.
\endproclaim 

As a function of $\mu =(-\lambda )^{\frac12}\in\Gamma $, $\frak A$ is
a strongly 
polyhomogeneous $\psi $do with symbol in $S^{0,0,1}(\Gamma )^{N\stimes
N}$ in local trivializations, its principal symbol being
equal to$$
\frak a^0(x',\xi ',\mu )=({p'}^0(x',\xi ')+\mu ^2)^{\frac12}.\tag2.3
$$
This follows essentially from Seeley \cite{S69}, since $p'(x',\xi ')+e^{2i\theta
}t^2$ is a classical elliptic symbol of order 2 with respect to the
cotangent variables $(\xi ',t)$. 
Moreover, $\frak A(-\mu ^2)$ is
invertible for $\mu \in\Gamma $, and parameter-elliptic (as defined
in \cite{G96}),  and 
$\frak A^{-1}$ has symbol in $S^{0,0,-1}(\Gamma )^{N\stimes N}$ in
local trivializations, with
principal part $(\frak a^0)^{-1}=({p'}^0(x',\xi ')+\mu ^2)^{-\frac12}$.
The symbols are holomorphic in $\mu
$. We observe furthermore that since $\partial _\lambda ^r\frak
A=c_r\frak A^{1-2r}$, $\partial _\lambda ^r\frak A$ has symbol in
$S^{0,0,1-2r}(\Gamma )^{N\stimes N}$ for $r\in\Bbb N$.

The indication by an upper index $N\stimes N$
means that the symbols are $N\stimes N$-matrix valued. 
The statement ``has symbol in
$S^{m,d,s}(\Gamma )^{N\times N}$ in local trivializations'' will be written briefly as:
``$\in \operatorname{OP}' S^{m,d,s}(\Gamma )$''.

In the case considered in Example 1.2, 
$\frak
A(\lambda )$ 
is the operator called $A_\lambda  $ in \cite{GS96} and \cite{G02}, $A_\mu $ in
\cite{GS95},
and $\frak A(0) =|A|$.

\proclaim{Definition 2.2}
The
Dirichlet-to-Neumann operator $A_{\operatorname{DN}}$ is defined
for $\lambda \notin \Sigma _{\delta ,R}$ by:$$
A_{\operatorname{DN}}(\lambda )=\gamma _1 K_{\operatorname{D}} (\lambda );\tag2.4
$$
cf\. {\rm (1.11)} ff.
\endproclaim 
 
The operator $A_{\operatorname{DN}}$ is a $\psi $do of order 1 for each
$\lambda $; this is a well-known fact in
the calculus of pseudodifferential boundary problems (cf\. Boutet de
Monvel \cite{BM71}, Grubb \cite{G96}). We observe moreover:
 
\proclaim{Lemma 2.3} The Dirichlet-to-Neumann operator
$A_{\operatorname{DN}}(\lambda )$ is a strongly polyhomogeneous $\psi 
$do in $E'$ of order $1$, which is principally equal to $-\frak A(\lambda )
$, i.e.,$$
A_{\operatorname{DN}}(\lambda )=-\frak A(\lambda )
+A'_{\operatorname{DN}}(\lambda ),\tag2.5
$$  
where $A_{\operatorname{DN}}(-\mu ^2)\in\operatorname{OP}' S^{0,0,1}(\Gamma )$ and
$A'_{\operatorname{DN}}(-\mu ^2)\in\operatorname{OP}'
S^{0,0,0}(\Gamma )$, with holomorphic symbols.
\endproclaim 

\demo{Proof} Clearly, $A_{\operatorname{DN}}(-\mu ^2)$ is strongly
polyhomogeneous of order 1, in the terminology of \cite{GS95}, \cite{G01},
since its symbol can be found from the corresponding 
calculation in the case where $\mu =e^{i\theta }t$ is replaced by $e^{i\theta }\partial
_{x_{n+1}}$. Then its symbol is in $S^{0,0,1}(\Gamma )^{N\stimes N}$.
Formula (2.5) is derived from the fact that the
principal part of 
$P+\mu ^2$ at $x_n=0$ equals the principal part of $-\partial
_{x_n}^2+P'+\mu ^2$.
For the latter operator considered on $X'\times \crp$,
the Poisson 
operator solving the semi-homogeneous Dirichlet problem as in (1.11) is
the mapping $\varphi (x')\mapsto z(x',x_n)=e^{-x_n\frak A(\lambda ) }\varphi
$ (as in \cite{G02, Prop\. 2.11} with $A^2$ replaced by $P'$);
application of $\partial _{x_n}$ followed by restriction to 
$x_n=0$ gives the mapping $\varphi \mapsto -\frak A(\lambda )\varphi
$. Then $A_{\operatorname{DN}}(\lambda )$ and $-\frak A(\lambda )$
are principally equal, so their difference $A'_{\operatorname{DN}}$
is strongly polyhomogeneous of order 0, hence has symbol in
$S^{0,0,0}(\Gamma )^{N\stimes N}$.
\qed \enddemo 

In particular, $A_{\operatorname{DN}}(-\mu ^2)$ is
parameter-elliptic, hence invertible for large enough $\mu $, the
inverse having symbol in $S^{0,0,-1}(\Gamma )^{N\stimes N}$ (see
Proposition 2.8 below). Note that $\partial _\lambda
^rA_{\operatorname{DN}}(\lambda )$ is strongly polyhomogeneous of
degree $1-2r$, hence lies in $\operatorname{OP}'S^{0,0,1-2r}(\Gamma )$.

In our construction of the resolvent below,
we need to be able to commute $A_{\operatorname{DN}}$
and $\Pi _1$ with an error having symbol in $S^{0,0,0}$. For this we
introduce

\proclaim{Assumption 2.4} The principal symbols of $\Pi _1$ and $P'$ commute.
\endproclaim 

This holds of course if $\Pi _1$ commutes with $P'$ (as in Example
1.2 with $P'=A^2$, $\Pi _1=\Pi _\ge$ or $\Pi _<$); it holds 
for general choices of $\Pi _1$ if $P'$ has scalar principal symbol.

\proclaim{Proposition 2.5} Under Assumption {\rm 2.4},
$[A_{\operatorname{DN}},\Pi _1]=A_{\operatorname{DN}}\Pi _1-\Pi _1\,
A_{\operatorname{DN}}$ and $[\frak A,\Pi _1]=\frak A\Pi _1-\Pi
_1\frak A$ are in $\operatorname{OP}' S^{0,0,0}(\Gamma )$, with
holomorphic symbols. 
Moreover, for any $r\ge 0$, the $r$'th
$\lambda $-derivatives are in $\operatorname{OP}' S^{0,0,-2r}(\Gamma )$. 
\endproclaim 

\demo{Proof} The main effort lies in the treatment of the case $r=0$.
Note first that 
since $A'_{\operatorname{DN}}$ and $\Pi _1$ are in
$\operatorname{OP}' S^{0,0,0}$, so is their commutator $[A'_{\operatorname{DN}},\Pi
_1]$, so in view of (2.5), what we have to show is that $[\frak A,\Pi
_1]$ is in
$\operatorname{OP}' S^{0,0,0}$. Denote $P'+\mu ^2={\overline P}$,
with symbol $p'(x ',\xi ')+\mu 
^2$ in local coordinates. The powers of ${\overline P}$ are defined for low
values of $s$ by$$
{\overline P}^s=\tfrac i{2\pi }\int_{\Cal C}\varrho  ^s({\overline P}-\varrho  )^{-1}\,d\varrho ,$$
where $\Cal C$ is a curve in $\Bbb C\setminus\crm$ encircling the
spectrum of  ${\overline P}$; we let it begin with a ray with angle $\delta $ and
end with a ray with angle $-\delta $, for some $\delta \in
\,]0,\frac\pi 2[\,$. (The location of $(r+\mu ^2)^{\frac12}$ when
$r\in \Bbb R_+$ is discussed in Remark 2.11 below.) Then since $[({\overline P}-\varrho )^{-1},\Pi _1]=({\overline P}-\varrho )^{-1}[\Pi
_1,P']({\overline P}-\varrho )^{-1}$, 
$$ 
[{\overline P}^s,\Pi _1]
=\tfrac i{2\pi }\int_{\Cal
C}\varrho  ^s({\overline P}-\varrho  )^{-1}[\Pi _1,P']({\overline P}-\varrho  )^{-1}\,d\varrho .
$$ 
Here, by the commutativity of the principal symbols,  $[\Pi _1,P']=M$ is
a first-order $\psi $do (independent of $\varrho $ and $\mu $). The integral
makes good sense for $s<1$ (converges in the norm of operators  
from $H^1(E')$ to $L_2(E')$), so we can write, since ${\overline P}
^\frac12=\frak A$,$$
[\frak A,\Pi _1]
=\tfrac i{2\pi }\int_{\Cal
C}\varrho  ^\frac12({\overline P}-\varrho  )^{-1}M({\overline P}-\varrho  )^{-1}\,d\varrho .\tag2.6
$$
The symbol of $[\frak A,\Pi _1]$ in a local coordinate system can be
found (modulo smoothing terms) from this formula. It is represented
by a series of terms obtained from (2.6) by insertion 
of the symbol expansions of the involved operators $({\overline P}-\varrho 
)^{-1}$ and $M$ and applications of the composition rule (2.2).
We have that the
symbol of $({\overline P}-\varrho  )^{-1}$  is an $(N\times N)$-matrix 
$r(x',\xi ',\mu ,\varrho )$ whose entries $r_{ij}$ are
series of strongly homogeneous functions  of $(\xi ',\mu
,\varrho ^\frac12)$, whereas the symbol $m(x',\xi ')$ of $M$ is a matrix
whose elements 
$m_{ij}$ are series of functions homogeneous in  $\xi '$ 
outside a neighborhhod of 0. The composition will give rise to integrals of
products 
$$\partial _{x',\xi '}^\gamma r_{ij}\,\partial _{x',\xi
'}^{\gamma '}m_{jk}\,\partial _{x',\xi '}^{\gamma ''}r_{kl};
$$
there are finitely many contributions to each degree of homogeneity.
The important thing is that in each homogeneous contribution, we can
take the factor coming from $m_{jk}$ outside the integral sign, since
it does not depend on $\varrho $ (nor on $\mu $). What is left is a
completely homogeneous integrand, which after integration gives a
strongly homogeneous function of $(\xi ',\mu )$. The factors coming
from $m$ are of degree $\le 1$, hence lie in $S^{1,0,0}$, and the
contributions from the integration 
are of degree $\le -1$, hence lie in $S^{0,0,-1}$ (moreover,  $\xi
'$-differentiation of order $\alpha $ lowers the former to $S^{1-|\alpha
|,0,0}$ and the latter to $S^{0,0,-1-|\alpha |}$). The full
contributions are then in $S^{1,0,0}\cdot S^{0,0,-1}\subset
S^{0,0,0}$, and can be collected in a series of terms of falling
degrees. Thus $[\frak A,\Pi _1]$ has a symbol series in $S^{0,0,0}$.

There is still the question of whether the remainder $\Cal R$, the
difference between the operator defined from $({\overline P}-\varrho  )^{-1}M({\overline P}-\varrho 
)^{-1}$ and an operator $S$ defined from a superposition of the homogeneous
terms, also has the right kind of symbol. It takes a certain effort to
prove this. We know \`a priori, since $\frak
A$ has symbol in $S^{0,0,1}$ and $\Pi _1$ has symbol in $S^{0,0,0}$,
that $[\frak A,\Pi _1]$ has symbol in $S^{0,0,1}\subset
S^{1,0,0}+S^{0,1,0}$. The remainder $\Cal R$ will be of order
$-\infty $ in this class, hence have symbol in $S^{-\infty ,1,0}$,
and we have to show that this can be reduced to $S^{-\infty ,0,0}$.

Recall that an operator $\Cal R$ with symbol in
$S^{-\infty ,1,0}$ has a kernel expansion $$
K(\Cal R)(x',y',\mu )\sim \sum _{j\ge
0}K_j(x',y')\mu ^{1-j},
$$ 
 with $K({\Cal R})-\sum_{j<J}K_j\mu ^{1-j}=O(\ang \mu
^{1-J})$ for all $J$,  where the $K_j(x',y')$ are $C^\infty $ kernels,
by \cite{GS95, Prop\. 1.21}. With $\Cal R_0$ denoting the operator with
kernel $K_0$, we can write $\Cal R=\mu \Cal R_0+\Cal R'$ where
$\Cal R'$ has the kernel expansion $\sum_{j\ge 1}K_j\mu ^{1-j}$,
hence symbol in $S^{-\infty ,0,0}$, so we are through if we show
that $\Cal R_0=0$. 

Let $\mu $ run on the ray
$\rp$, and let $\varepsilon \in \,]0,\frac12[\,$. Then
$$
{\overline P}^{-\varepsilon }[\frak A,\Pi _1]
=\tfrac i{2\pi }\int_{\Cal
C}\varrho  ^\frac12({\overline P}-\varrho  )^{-1}{\overline P}^{-\varepsilon }M({\overline P}-\varrho  )^{-1}\,d\varrho .
$$
Since ${\overline P}=P'+\mu ^2$,  $\|{\overline P}^sf\|_{L_2}\simeq \|f\|_{H^{2s,\mu }}$
(uniformly in $\mu $ for $\mu \ge 1$), where $H^{t,\mu }$ is the
Sobolev space defined by localization from
the case $H^{t,\mu }(\Bbb R^{n-1})$ with norm $\|u\|_{H^{t,\mu
}}=\|\ang{(\xi ',\mu )}^t\hat u(\xi ')\|_{L_2}$, $t\in\Bbb R$. 
From the resolvent
estimates $$
\|\varrho ({\overline P}-\varrho )^{-1}f\|_{L_2}\le C_1\|f\|_{L_2},\quad \|({\overline P}-\varrho )^{-1}\|_{H^{2,\mu
}}\simeq \|{\overline P}({\overline P}-\varrho )^{-1} f\|_{L_2}\le C_2\|f\|_{L_2},\tag2.7
$$
that are valid for $\varrho \in\Cal C$,
then follow by interpolation  that $$
|\varrho |^{\frac12+\varepsilon }\|({\overline P}-\varrho )^{-1}f\|_{H^{1-2\varepsilon ,\mu }
}\le C_3\|f\|_{L_2}.\tag2.8
$$
We also have that $$
\|{\overline P}^{-\varepsilon }Mv\|_{L_2}\le C_4\|Mv\|_{H^{-2\varepsilon ,\mu
}}\le 
C_5\|Mv\|_{H^{-2\varepsilon }}\le C_6\|v\|_{H^{1-2\varepsilon }}\le
C_7\|v\|_{H^{1-2\varepsilon ,\mu }},\tag 2.9
$$ using that $1-2\varepsilon \ge 0$.
Thus
$$\multline
\|{\overline P}^{-\varepsilon }[\frak A,\Pi _1]\|_{\Cal L(L_2)}\\
\le
C\int_{\Cal
C}|\varrho  |^\frac12\|({\overline P}-\varrho  )^{-1}\|_{\Cal L(L_2)}\|{\overline P}^{-\varepsilon
}M\|_{\Cal L(H^{1-2\varepsilon,\mu } ,L_2)}\|\|({\overline P}-\varrho )^{-1}\|_{\Cal 
L(L_2,H^{1-2\varepsilon ,\mu })}\,|d\varrho |\\
\le C'\int_{\Cal C}|\varrho |^{\frac12-1-\frac12-\varepsilon }\,|d\varrho |\le C'',
\endmultline\tag2.10$$ 
where we applied (2.7) to the first factor, (2.9) to the middle
factor and (2.8) to the last factor.
We have that $[\frak A,\Pi _1]=S+\mu \Cal R_0+\Cal R'$, where $S+\Cal R'$
is bounded 
in $L_2$ uniformly in $\mu $ for $\mu \ge 1$, hence ${\overline P}^{-\varepsilon }(S+\Cal R')$ is \`a
fortiori so. Thus, by (2.10), $\|\mu \Cal R_0f\|_{H^{-2\varepsilon ,\mu
}}$ is bounded in $\mu $ for $\mu \ge 1$, for any $f$. But if $\Cal R_0f$
is a nonzero function, $\mu \|\Cal R_0f\|_{H^{-2\varepsilon ,\mu }}$
cannot be bounded for $\mu \to \infty $ (recall that $\varepsilon
<\frac12$). Thus $\Cal R_0=0$, and the proof for the case $r=0$ is
complete.

For $r\ge 1$, we observe that $\partial _\lambda ^r\frak A=c_r\frak
A^{1-2r}$ and that 
$$\aligned
[\Pi _1,\frak A^{-1}]&=\frak A^{-1}\frak A\Pi _1\frak A^{-1}-
\frak A^{-1}\Pi _1\frak A\frak A^{-1}=\frak A^{-1}[\frak A,\Pi
_1]\frak A^{-1}\in \operatorname{OP}' S^{0,0,-2};\\
[\Pi _1,\frak A^{-k}]&=\Pi _1\frak A^{-1}\frak A^{-k+1}-\frak
A^{-1}\Pi _1\frak A^{-k+1}+\frak
A^{-1}\Pi _1\frak A^{-k+1}-\dots 
-\frak A^{-k}
\Pi _1\\
&=\sum_{0\le j\le k-1}\frak
A^{-j}[\Pi _1,\frak A^{-1}]\frak A^{-k+1+j}\in \operatorname{OP}' S^{0,0,-1-k}
\endaligned\tag 2.11$$
by the result for $r=0$ and the composition rules. The first
calculation shows the assertion 
on $\partial _\lambda ^r[\Pi _1,\frak A]$
for $r=1$; the second calculation shows it for general $r\ge 2$, when
we take $k=1-2r$. 

The result now likewise follows for $\partial _\lambda ^r[\Pi
_1,A_{\operatorname{DN}}]$, when we
use that $\partial _\lambda ^rA_{\operatorname{DN}}=-\partial
_\lambda ^r\frak A+\partial _\lambda ^rA'_{\operatorname{DN}}$, where
$\partial _\lambda ^rA'_{\operatorname{DN}}\in \operatorname{OP}' S^{0,0,-2r}$ since it
is strongly polyhomogeneous of degree $-2r$.
\qed
\enddemo

The resolvent $(P_T-\lambda )^{-1}$ is the solution operator for the problem
$$\aligned
(P-\lambda )u&=f\text{ on }X,\\
\Pi _1\gamma _0u&=0\text{ on }X',\\
\Pi_2(\gamma _1u+B\gamma _0u)&=0\text{ on }X',
\endaligned\tag2.12$$
where $\lambda $ runs in a suitable subset of $\Bbb C$.
The problem of constructing the resolvent will be transformed by some
auxiliary constructions:  

First, 
we can extend $P$ to a strongly elliptic principally selfadjoint
differential operator $\widetilde P$ 
on an $n$-dimensional compact manifold $\widetilde X$ in which $X$ is smoothly
imbedded, modifying the definition of $R(\delta )$ in (1.9) such that 
also the spectrum of $\widetilde P$ is contained in $\Sigma _{\delta
,R}$. Assume in the following that $\lambda \notin \Sigma _{\delta
,R}$. Let $Q(\lambda )$ 
be the inverse of $\widetilde P-\lambda $ and let $Q_+=r^+Qe^+$ be
its truncation to $X$ ($r^+$ is restriction from $\widetilde X$ to
$X$, $e^+$ is extension by zero on $\widetilde X\setminus X$). 
It is well-known that 
$R_{\operatorname{D}}(\lambda )$ and $K_{\operatorname{D}}(\lambda )$ (cf.\
(1.10), (1.11)) are connected by the formula
$$
R_{\operatorname{D}}=Q_{+}-K_{\operatorname{D}}\gamma _0Q_+,\tag 2.13
$$
for $\lambda \notin \Sigma _{\delta ,R}$.
Let $$ 
v=R_{\operatorname{D}}f,\quad z=u-v,\quad  \psi =-\Pi _2\gamma _1v,\tag2.14
$$ 
then (2.12) may be replaced by the problem for $z$:
$$\aligned
(P-\lambda )z&=0\text{ on }X,\\
\Pi _1\gamma _0z&=0\text{ on }X',\\
\Pi_2(\gamma _1z+B\gamma _0z)&=\psi \text{ on }X'.
\endaligned\tag2.15$$
Here one can let $f\in L_2(E)$ so that the solution $u$ is sought in
$H^2(E)$; then $v\in H^2(E)$, $z\in H^2(E)$, and $\gamma _jv$ and $\gamma
_jz$ are in $H^{\frac32-j}(E')$, and so are $\Pi _i\gamma _jv$ and
$\Pi _i\gamma _jz$, since the projections $\Pi _j$ are bounded in $H^s(E')$
for all $s\in\Bbb R$. One could also carry out the whole calculation for
$C^\infty $ sections while keeping track of the orders of the
involved operators. 

We now make a reduction of (2.15) using $A_{\operatorname{DN}}$. When $z\in H^2(E)$
with $(P-\lambda )z=0$,
$z$ is in 1-1 correspondence with $\gamma _0z$ running through the
space $H^{\frac32}(E')$:
$$
z=K_{\operatorname{D}} \gamma _0z,\tag2.16
$$
in view of the unique solvability of the Dirichlet problem.
Denote $\gamma _0z=\varphi $; then we can replace $\gamma _0z$ by
$\varphi $ in (2.15). Moreover, $\gamma 
_1z=A_{\operatorname{DN}} \varphi $. So, with $z=K_{\operatorname{D}}
\varphi $ 
(assuring the validity of the
first line in (2.15)), we arrive at the problem for $\varphi $:
$$\aligned
\Pi _1\varphi &=0,\\ 
\Pi_2(A_{\operatorname{DN}} 
+B)\varphi 
&=\psi . 
\endaligned \tag2.17
$$ 

Since the mapping $\gamma _1\: D(P_{\operatorname{D}})\to
H^{\frac12}(E')$ is surjective, $\gamma _1v$ runs through all of
$H^{\frac12}(E')$ when $f$ runs through $L_2(E)$. So the reduced
problem (2.17) must be solved for all $\psi \in \Pi
_2(H^{\frac12}(E'))$. 
When $f$ is given in $L_2(E)$ and $v$ and $\psi $
are defined by (2.14), then, if $\varphi \in H^{\frac32}(E')$ solves
(2.17) and we set $z=K_{\operatorname{D}}\varphi$, $  u=v+z$, we find
that $u\in H^2(E)$ solves (2.12).

\proclaim{Lemma 2.6} For $\lambda \notin\Sigma _{\delta ,R}$, define
$$\aligned 
S(\lambda )&=A_{\operatorname{DN}} +[A_{\operatorname{DN}},\Pi _1]
+\Pi_2 B\Pi_2,\\
S'(\lambda )&=A_{\operatorname{DN}} -[A_{\operatorname{DN}},\Pi _1]
.\endaligned\tag 2.18$$
When $\psi $ runs through $\Pi _2H^{\frac12}(E')$, the problem {\rm
(2.17)} is uniquely solvable with $\varphi \in H^{\frac32}(E')$ for
those $\lambda  
$ for which $S(\lambda )$ and $S'(\lambda )$ are invertible from
$H^{\frac32}(E')$ to $H^{\frac12}(E')$.
When this holds, the solution is$$
\varphi =S(\lambda )^{-1}\psi  .\tag2.19$$
Uniqueness of solution holds when merely $S(\lambda )$ is injective.
\endproclaim 

\demo{Proof}
Assume first that $\varphi $ is a solution of (2.17).
By the first equation, $\varphi =\Pi_2 \varphi $, so we
can rewrite the second equation
as:
$$
A_{\operatorname{DN}}\Pi_2\varphi
+[\Pi_2,A_{\operatorname{DN}}]\varphi +\Pi_2 B\Pi_2\varphi =\psi
 .\tag2.20
$$
Now we compose the first equation of (2.17) with
$A_{\operatorname{DN}} $ and add it to 
(2.20); this gives that $\varphi $ satisfies
$$
(A_{\operatorname{DN}}+[\Pi _2,A_{\operatorname{DN}}]+\Pi _2 B\Pi _2)\varphi
=\psi .
$$
Since
$[\Pi_2,A_{\operatorname{DN}}]=[A_{\operatorname{DN}},\Pi _1]$, the
operator in the left-hand side equals $S(\lambda )$ defined in
(2.18), so $\varphi $ satisfies
$$
S(\lambda )\varphi =\psi .\tag2.21
$$
When $S(\lambda )$ is injective, there is at most one solution
$\varphi $ of (2.21), hence of (2.17); when $S(\lambda )$ is
bijective, the only possibility is that it equals 
$S(\lambda )^{-1}\psi $. 

Now assume that also $S'(\lambda )$ is invertible, and, for a given
$\psi $ with $\Pi _1\psi =0$, define $\varphi $ by (2.19). Then of
course (2.21) holds, and we shall show that $\varphi $ solves (2.17).
Since $\Pi _1\Pi _2=0$
and $\Pi _1^2=\Pi _1$, an application of $\Pi _1$ to (2.21) gives:
$$\aligned
0&=\Pi _1S\varphi =\Pi _1(A_{\operatorname{DN}}\varphi
+[A_{\operatorname{DN}},\Pi _1]\varphi +\Pi_2 B\Pi_2\varphi )
\\&=\Pi _1A_{\operatorname{DN}}\varphi +\Pi _1A_{\operatorname{DN}}\Pi
_1\varphi -\Pi _1^2A_{\operatorname{DN}}\varphi =
\Pi _1A_{\operatorname{DN}}\Pi
_1\varphi\\
&=(A_{\operatorname{DN}}-[A_{\operatorname{DN}},\Pi _1])\Pi _1\varphi
=S'(\lambda )\Pi _1\varphi .\endaligned\tag2.22$$
It follows that $\Pi _1\varphi =0$, so
$\varphi $ satisfies the first equation in (2.17). The second
equation is then retrieved from this and (2.21): Since $\Pi _2\varphi
=\varphi $,
$$\aligned
\psi &=(A_{\operatorname{DN}} +[A_{\operatorname{DN}},\Pi _1]
+\Pi_2 B\Pi_2)\varphi \\
&=A_{\operatorname{DN}}\Pi _2\varphi +[\Pi
_2,A_{\operatorname{DN}}]\varphi 
+\Pi_2 B\varphi =\Pi _2(A_{\operatorname{DN}} +B)\varphi .\quad\square
\endaligned\tag2.23$$
\enddemo 

When Assumption 2.4 holds, the commutator 
terms in (2.18)  have symbol in $S^{0,0,0}$. Then we will show that
$S'(\lambda )$ can be inverted for large $|\lambda
|$ within the weakly 
polyhomogeneous calculus, and that the same holds for $S(\lambda )$ when 
the following assumption is satisfied:

\proclaim{Assumption 2.7} There is a $\theta \in \,]0,\frac\pi 2]$ such
that, with  $b ^h(x',\xi  ')$ and  $\pi _i ^h (x',\xi 
')$  denoting the strictly homogeneous principal
symbols of 
$B$ and the $\Pi _i$,$$
\frak a^0 -\pi _2^h b^h  \pi ^h
_2\text{ is
invertible  for }\xi '\in\Bbb R^n,\mu \in \Gamma
_{\theta}\cup\{0\}\text{ with }(\xi ',\mu
)\ne (0,0).\tag2.24
$$ 

\endproclaim  

This will allow construction of a resolvent family for large $\lambda
$ with $|\arg(-\lambda )|<2\theta $. If, moreover,
$\theta >\frac\pi 4 $, we can also construct a heat operator family.
Note that when  $B$ is a differential operator, its
symbol is a polynomial, so $b^h$ equals the usual principal symbol $b^0$. For
$\Pi _2$, the strictly homogeneous 
principal symbol $\pi _2^h$ is in general not continuous at $\xi '=0$, but when it is
multiplied by $b^h$, which is $O(|\xi '|)$, we get a continuous
function at $\xi '=0$ (taking the value 0 there).

\proclaim{Proposition 2.8} Let Assumption
{\rm 2.4} hold.

{\rm (i)} For each $\theta '\in \,]0,\frac\pi 2[\,$ there is an
$r(\theta ')\ge 0$ such that
$A_{\operatorname{DN}}(-\mu ^2)$ and $S'(-\mu ^2 )$ are invertible for
$\mu \in \Gamma _{\theta ',r}$.

{\rm (ii)} Let moreover Assumption {\rm 2.7} hold. Then 
for
each $\theta '\in \,]0 ,\theta [\,$ there is an $r(\theta ')\ge 0$ such
that $S(-\mu ^2)$ is invertible for $\mu \in \Gamma _{\theta ',r}$.
 
Here the operator families
$A_{\operatorname{DN}}(-\mu ^2 )^{-1}$ and $S'(-\mu ^2 )^{-1}$ belong
to
$\operatorname{OP}'S^{0,0,-1}(\Gamma )$, and 
$S(-\mu ^2)^{-1}$ belongs to $\operatorname{OP}'S^{0,0,-1}(\Gamma _\theta
)$, with holomorphic symbols. For each $r\in\Bbb N$, $\partial
_\lambda ^r$ map them into operators in
$\operatorname{OP}'S^{0,0,-1-2r}(\Gamma _\theta )$.

\endproclaim 

\demo{Proof}
Let us go directly to the proof of (ii), the statements in (i) are proved
by easier variants.
We have on one hand that $S(-\mu ^2)$ is composed
of operators of the type considered in \cite{G96}, with  the $\mu 
$-dependent factors of 
regularity $\infty $ and the $\Pi _i$ of regularity $0$, $B$ of
regularity $\infty $ resp\. 1 if it is a differential resp\.
pseudodifferential operator. Removing the commutator term
$[A_{\operatorname{DN}},\Pi _1]$ 
from (2.20) for a 
moment, we have an operator$$
S''(-\mu ^2 )=A_{\operatorname{DN}}
+\Pi_2 B\Pi_2=-\frak A+A'_{\operatorname{DN}}+\Pi_2 B\Pi_2
, \tag2.25
$$
which is of regularity 1 (since $\Pi_2 B\Pi_2$ is so), so that the
invertibility of the  strictly homogeneous principal symbol $-\frak
a^0+\pi ^h_2b^h\pi _2^h $
assures parameter-ellipticity in the sense of \cite{G96}, cf\. Prop\.
2.1.12 there. Then $S''(-\mu ^2 )$ is
invertible for $ \mu\in \Gamma _{\theta ',r}$  with a sufficiently
large $r$, the inverse being continuous 
from $H^{s,\mu }(E)$ to $H^{s+1,\mu }(E)$ for $s\in \Bbb R$, by \cite{G96,
Th\. 3.2.11}. Since the commutator 
term has $L_2$-norm bounded in $\mu  $ by Proposition 2.5, we
find by a Neumann 
series argument that $S(-\mu ^2 )$ itself is invertible for large $r$
with an inverse that is bounded from $L_2(E)$ to $H^{1,\mu }(E)$.

Now $\frak A^{-1}S$ is likewise invertible for
the considered $\mu $, and lies in $\operatorname{OP}'S^{0,0,0}(\Gamma _\theta
)$ since $\frak A^{-1}$ lies in $\operatorname{OP}'S^{0,0,-1}(\Gamma 
)$. Then the ``spectral invariance theorem'' \cite{G99, Th\.
6.5}, applied to $\frak 
A^{-1}S$, shows that its 
inverse $S^{-1}\frak A$ belongs to our weakly polyhomogeneous
calculus and lies in 
$\operatorname{OP}'S^{0,0,0}(\Gamma _\theta )$. It follows that
$S^{-1}=(S^{-1}\frak A)\frak A^{-1} \in 
\operatorname{OP}'S^{0,0,-1}(\Gamma _\theta )$. Since holomorphy is
preserved under composition, the resulting symbols are holomorphic.
The statements on $\lambda $-derivatives follows by successive applications of
the formula $\partial _\lambda S^{-1}=-S^{-1}\partial _\lambda S\,
S^{-1}$, using that $\partial _\lambda S=\partial _\lambda
A_{\operatorname{DN}}+\partial _\lambda [A_{\operatorname{DN}},\Pi
_1]$, with properties described in Lemma 2.3 ff.\ and Proposition 2.4.

Similar proofs work for the other operators without the complication
due to the presence of $B$  (so any $\theta '\in \,]0,\frac\pi 2[\,$ is allowed there).
\qed
\enddemo 

\example{Remark 2.9}
It should be noted that the proof of Lemma 2.6 only shows the
necessity of unique 
solvability of (2.21) for $\psi $ in the range of $\Pi _2$, so that
Assumption 2.7 (assuring solvability for general $\psi $) may
seem too strong. However, when Assumption 2.4 holds, solvability of
(2.21) for $\psi \in \Pi _2H^{\frac12}(E')$, $\mu =re^{i\theta _0}$,
$r\ge r_0$ (some $r_0\ge 0$, $|\theta _0|<\theta $), implies
solvability for all $\psi \in H^{\frac12}(E')$, $r\ge r_1$ with some
$r_1\ge r_0$:

Uniqueness of course holds. Existence is assured as follows: Let
$\psi \in H^{\frac12}(E')$ and write $\psi =\Pi _1\psi +\Pi _2\psi $.
Define$$
\widetilde S(\lambda )\psi =A_{\operatorname{DN}}^{-1}\Pi _1\psi +S^{-1}\Pi
_2\psi ,
$$
for $r$ so large that also $A_{\operatorname{DN}}$ is invertible.
Then$$\aligned
S\widetilde S\psi &=SA_{\operatorname{DN}}^{-1}\Pi _1\psi +\Pi
_2\psi 
=\psi +[A_{\operatorname{DN}},\Pi _1]A_{\operatorname{DN}}^{-1}\Pi
_1\psi +\Pi _2B\Pi _2A_{\operatorname{DN}}^{-1}\Pi _1\psi \\
&=\psi +[A_{\operatorname{DN}},\Pi _1]A_{\operatorname{DN}}^{-1}\Pi
_1\psi +\Pi _2B\Pi _2[A_{\operatorname{DN}}^{-1},\Pi _1]\psi 
=(I+\widetilde S_1)\psi .
\endaligned$$
where $\widetilde S_1=[A_{\operatorname{DN}},\Pi _1]A_{\operatorname{DN}}^{-1}\Pi
_1 +\Pi _2B\Pi _2A_{\operatorname{DN}}^{-1}[\Pi
_1,A_{\operatorname{DN}}]A_{\operatorname{DN}}^{-1}$ has symbol in
$S^{0,0,-1}(\Gamma )$. So, like the operators treated in the proof of
Proposition 2.8, $I+\widetilde S_1$ is invertible for large enough $\mu $ on the
ray, and $\varphi =\widetilde S(I+\widetilde S_1)^{-1}\psi $ solves the equation $S\varphi =\psi $.

\endexample

We can now describe the resolvent. Here we shall use the terminology of weakly
polyhomogeneous pseudodifferential boundary operators worked out in
\cite{G01} (the relevant parts summed up in \cite{G02, Sect\. 2}), extending the
calculus of Boutet de Monvel \cite{BM71}. One can get
quite far with linear combinations of compositions of elementary
operators as in \cite{GS95}, \cite{G99}, but when the expressions get
increasingly complicated, it seem advantageous to use the systematic
calculus.
We shall here just recall the basic definitions in the case where
$X$, $X'$ are replaced by $\overline{\Bbb R}^n_+$, $\Bbb R^{n-1}$.

One considers Poisson operators $K$ (mapping functions on $\Bbb R^n_+$ to
functions on $\Bbb R^{n-1}$), trace operators $T$ of class 0 (mapping
functions on  
$\Bbb R^{n-1}$ to functions on $\Bbb R^n_+$) and singular Green operators
$G$ of class 0 (mapping functions on  $\Bbb R^n_+$
to functions on $\Bbb R^n_+$), of the form
$$\aligned
K&=\operatorname{OPK}(\tilde k)\:v(x')\mapsto \int_{\Bbb
R^{2(n-1)}}e^{i(x'-y')\cdot
\xi '} \tilde 
k(x',x_n,\xi ',\mu )v(y')\,dy'\d\xi ',\\
T&=\operatorname{OPT}(\tilde t)\:u(x)\mapsto \int_{\Bbb
R^{2(n-1)}}\int_0^\infty e^{i(x'-y')\cdot
\xi '} \tilde 
t(x',x_n,\xi ',\mu )u(y',x_n)\,dx_ndy'\d\xi ',\\
G&=\operatorname{OPG}(\tilde g)\:u(x)\mapsto \int_{\Bbb
R^{2(n-1)}}\int_0^\infty e^{i(x'-y')\cdot
\xi '} \tilde 
g(x',x_n,y_n,\xi ',\mu )u(y)\,dy\d\xi '.
\endaligned\tag2.26$$
Note that the usual $\psi $do definition is used with respect to the
$x'$-variable. In fact, we can view OPK, OPT and OPG as
$\operatorname{OPK}_n\operatorname{OP}'$,
$\operatorname{OPT}_n\operatorname{OP}'$  resp\.
$\operatorname{OPG}_n\operatorname{OP}'$, where
$\operatorname{OPK}_n$ etc.\ stand for the application of (2.26)
with respect to $x_n$-variables alone.

The functions $\tilde k$, $\tilde t$ and $\tilde g$ are called the
{\it symbol-kernels} of $K$, $T$, resp\. $G$. There is also a ``complex
formulation'', where e.g.\ the symbol-kernel $\tilde k(x',x_n,\xi
',\mu )$ is replaced by the {\it symbol} $k(x',\xi ',\xi _n,\mu
)=\Cal F_{x_n\to \xi 
_n}\tilde k$, and $\Cal F^{-1}_{\xi _n\to x_n}$ is included in the
definition of the operator ($\Cal F$ denotes Fourier transformation).

We say that $\tilde k\in\Cal
S^{m,d,s}(\Lambda ,\Cal S_+)$, resp\. $\tilde g \in\Cal
S^{m,d,s}(\Lambda ,\Cal S_{++})$, when 
$$\gather
\sup_{u_n\in \Bbb R_+}| \partial 
_{|z|}^j(z^{d}\kappa ^{-s-1}\ang{z\xi '}^{l-l'}u_n^l\partial _{u_n}^{l'}\tilde
k({x'},|z|u_n,{\xi '},1/z)) |
 \le C\ang{\xi 
'}^{m+j},\text{ resp.} \tag2.27\\
\sup_{u_n,v_n\in \Bbb R_+}| \partial 
_{|z|}^j(z^{d}\kappa ^{-s-2}\ang{z\xi '}^{l-l'+k-k'}u_n^l\partial
_{u_n}^{l'}
v_n^k\partial _{v_n}^{k'}\tilde
g({x'},|z|u_n,|z|v_n,{\xi '},1/z)) |
 \le C\ang{\xi 
'}^{m+j} ,
\endgather
$$ 
for all indices,
uniformly for $|z|\le 1$, $\Zfrac$ in closed subsectors of $\Lambda $,
with similar estimates for the derivatives $\partial _{\xi '}^\alpha
\partial _{x'}^\beta $ with $m$ replaced by $m-|\alpha |$. Here
$\kappa =(|\mu |^2+|\xi '|^2)^{\frac12}= (|\Zfrac|^2+|\xi '|^2)^{\frac12}$.  
 $\tilde t$ is estimated like $\tilde k$. 

Again the symbol-kernels are said to be holomorphic (in $\mu $), when
they and their derivatives are holomorphic for $\mu \in\Lambda
^\circ$, $|\mu |+|\xi '|\ge \varepsilon >0$. The motivation for the
scaling $x_n=|z|u_n$ is explained in \cite{G01}, where complete and
satisfactory composition rules are worked out.
The operators are defined on $X$, $X'$ by standard localization
methods. 

We recall one
further operation, that of taking the normal trace $\tr_n$: When $G$ is a
singular Green operator as above, the normal trace $\tr_n G$ is the
$\psi $do on $\Bbb R^{n-1}$ with symbol
$$
(\operatorname{tr}_n\tilde g)(x',\xi ',\mu )=\int_0^\infty \tilde g(x',x_n,x_n,{\xi '},\mu
)\,dx_n\tag 2.28$$
(called $\ttilde g$ in \cite{G96}). Here the symbol map $\tr_n$ acts as follows:$$
\tr_n\: \Cal
S^{m,d,s-1 }(\Lambda  ,\srpp)\to S^{m,d,s}(\Lambda ).\tag2.29
$$
For operators of trace class,
$$
\Tr_{\rnp}G=\Tr_{\Bbb R^{n-1}}(\tr_n G)\tag 2.30
$$
(if $G$ has the kernel $K(x,y,\mu )$ then $\tr_nG$ has the kernel
$\int_0^\infty K(x',x_n,y',x_n,\mu )\,dx_n$),
and there is a similar rule for the operators carried over to the
manifold situation, when the symbol-kernel of $G$ is supported in
$X_c$ and the volume element on $X_c$ is taken of the form $v(x')dx'dx_n$:$$
\Tr_XG=\Tr_{X'}(\tr_nG).\tag2.31
$$
One has that  $G-\chi G\chi $ is smoothing and
$O(|\mu | ^{-M})$ for $\mu \to\infty $ in closed subsectors of
$\Lambda $, all $M$,  
when $\chi \in C_0^\infty (\,]-c,c[\,)$ and is 1 near
$x_n=0$ (cf\. e.g.\ \cite{G01, Lemma 7.1}). Thus in the trace expansion
calculations for singular Green operators on $X$, we can replace $G$
by $\chi G\chi $ and use (2.31) to reduce to a calculation for a
$\psi $do on $X'$; here \cite{GS95, Th\. 2.1} can be applied.

We shall express the fact that ``the
symbol-kernel is in $\Cal S^{m,d,s}(\Gamma _\theta ,\Cal
S_{+})^{N\times N}$
resp\. \linebreak$\Cal S^{m,d,s}(\Gamma _\theta ,\Cal S_{++})^{N\times N}$ in local
trivializations'' more briefly by saying that the operator lies in
$\operatorname{OPK}\Cal S^{m,d,s}(\Gamma _{\theta } ,\Cal S_{+})$,
$\operatorname{OPT}\Cal S^{m,d,s}(\Gamma _{\theta } ,\Cal S_{+})$
resp\. $\operatorname{OPG}\Cal S^{m,d,s}(\Gamma _{\theta } ,\Cal
S_{++})$.

The symbol-kernels we consider, moreover have expansions in
appropriately quasi-homogeneous terms, e.g., $\tilde
g\sim\sum_{j\in\Bbb N}\tilde 
g_{m-j}$ with  $\tilde g_{m-j}\in \Cal S^{m-j,d,s}(\Gamma _{\theta } ,\Cal
S_{++})$. There is the same distinction between weakly
polyhomogeneous symbol-kernels
(with
homogeneity 
 for $|\xi '|\ge 1$)
and strongly  polyhomogeneous symbol-kernels
(with homogeneity for $|\xi '|+|\mu |\ge 1$, etc\.) as for $\psi $do symbols.

The following elementary examples are basic in our calculations:
\roster
\item
$K_D(-\mu ^2)$ is a strongly polyhomogeneous Poisson operator 
in $\operatorname{OPK}\Cal S^{0,0,-1}(\Gamma ,\Cal
S_{+})$.
\item
$\gamma _0Q(-\mu ^2)_+$ and $\gamma _1Q(-\mu ^2)_+$ are strongly
polyhomogeneous trace 
operators of class 
0 
 in $\operatorname{OPT}\Cal S^{0,0,-2}(\Gamma ,\Cal
S_{+})$ resp\. $\operatorname{OPT}\Cal S^{0,0,-1}(\Gamma ,\Cal
S_{+})$.
\endroster
As mentioned in the proof of Lemma 2.3, $K_{\operatorname{D}}$ is
principally the same as the operator $K_{\frak A}\: v\mapsto
e^{-x_n\frak A}v$, 
for $x_n\in X_c$ (cf\. \cite{G02, Prop\. 2.11}).
We also have (by \cite{G02, (1.17), (4.14)} with $A^2$ replaced by $P'$)
that $\gamma _0Q_+$ and $\gamma 
_1Q_+$ {\it act principally,} for functions supported in $X_c$, {\it like the
operators $\frac12 \frak A^{-1}T_{\frak A}$ resp\.} $\frac12 T_{\frak
A}$, where  $$
T_{\frak A}u=\int_0^\infty e^{-x_n\frak A}u(x',x_n)\,dx_n.\tag2.32
$$ 
Simple examples of 
singular Green operators are $K_{\operatorname{D}}\gamma _0Q_+$ (whose
negtive is the singular Green part of
$R_{\operatorname{D}}$, cf\. (2.13)) and $K_{\operatorname{D}}\gamma
_1Q_+$; they are strongly polyhomogeneous
and belong to $\operatorname{OPG}\Cal S^{0,0,-3}(\Gamma ,\Cal
S_{++})$ resp\. $\operatorname{OPG}\Cal S^{0,0,-2}(\Gamma ,\Cal
S_{++})$. $2K_{\operatorname{D}}\gamma _1Q_+$  is (on $X_c$) principally equal to $G_{\frak A}$, which
acts as follows:
$$
G_{\frak A}u=\int_0^\infty e^{-(x_n+y_n)\frak A}u(x',y_n)\,dy_n.\tag2.33 
$$
For the latter, $\tr_n$ is easy to determine by functional calculus:
$$
\tr_n G_{\frak A}=\int_0^\infty e^{-2x_n\frak A}\,dx_n=(2\frak
A)^{-1};\tag2.34  
$$
this kind of calculation plays a role in our analysis of trace coefficients in
Section 4.

\proclaim{Theorem 2.10} Let Assumptions {\rm 1.1}, {\rm 2.4} and {\rm
2.7} hold.  Then for
each $\theta '\in \,]0,\theta [\,$ there is an $r=r(\theta ')\ge 0$ such
that for $\mu \in \Gamma _{\theta ', r}$,
$P_T+\mu ^2=P_T-\lambda $ is a bijection from 
$D(P_T)$ to $L_2(E)$ with inverse $(P_T-\lambda )^{-1}=R_T(\lambda )$
of the form $$
\aligned
R_T(\lambda )&=Q(\lambda )_++G(\lambda ),\\
G(\lambda )&=-K_{\operatorname{D}}(\lambda )\gamma _0Q(\lambda )_++K_{\operatorname{D}}(\lambda ) [S_0(\lambda )\gamma
_0+S_1(\lambda )\gamma _1]Q(\lambda )_+;
\endaligned\tag2.35
$$ 
here $S_0$ and $S_1$ (given in {\rm (2.36)} below) are weakly
polyhomogeneous $\psi $do's in $E'$ 
lying in $\operatorname{OP}'S^{0,0,0}(\Gamma _{\theta })$ 
resp\. $\operatorname{OP}'S^{0,0,-1}(\Gamma _{\theta })$,  and hence
$G$ is a singular Green operator of class $0$ 
in $\operatorname{OPG}\Cal
S^{0,0,-3}(\Gamma _{\theta },\Cal 
S_{++})$.
\endproclaim 

\demo{Proof} 
For a $\theta '\in\,]0,\theta [\,$, take $r=r(\theta ')$ such that
$A_{\operatorname{DN}}$ is well-defined and the operators $S(\lambda
)$ and $S'(\lambda )$ are invertible for $\lambda =-\mu ^2$, $\mu
\in\Gamma _{\theta ',r}$.
Then (2.17) is solved uniquely by (2.19), and it remains to draw the
conclusions for the original problem (2.12). In view of (2.13),
(2.14) and (2.19), the solution is:
$$\aligned
u&=v+z=R_{\operatorname{D}}f+K_{\operatorname{D}} \varphi \\
&=R_{\operatorname{D}}f-K_{\operatorname{D}} S^{-1}\Pi _2\gamma _1R_{\operatorname{D}}f\\
&=Q_+f-K_{\operatorname{D}}\gamma _0Q_+f- K_{\operatorname{D}}
S^{-1}\Pi _2 \gamma _1 Q_+f +K_{\operatorname{D}}
S^{-1}\Pi _2 A_{\operatorname{DN}}\gamma _0 Q_+f.
\endaligned$$
This shows (2.35) with$$
S_0=S^{-1}\Pi _2A_{\operatorname{DN}},\quad S_1=-S^{-1}\Pi _2;\tag 2.36
$$ 
they lie in $\operatorname{OP}'S^{0,0,0}(\Gamma _\theta )$ resp\.
$\operatorname{OP}'S^{0,0,-1}(\Gamma _\theta )$  by the rules of
calculus. The statement on $G$ now follows from the information given
before the theorem on $K_{\operatorname{D}}$, $ \gamma _0Q_+$ and 
$ \gamma _1Q_+$, and the composition rules.
\qed
\enddemo

\example{Remark 2.11} Let us give some sufficient conditions for the
validity of Assumption 2.7. Consider, in a local trivialization, a
point $(x',\xi ')$ with $|\xi 
'|=1$ (the result is carried over to general $\xi '\ne 0$ by
homogeneity, and for $\xi '=0$ the assumption is trivially
satisfied). Let $$ 
0<\lambda _1(x',\xi ')\le \lambda _2(x',\xi ')\le\cdots\le \lambda
_N(x',\xi ')  
$$
be the eigenvalues of the matrix ${p'}^0(x',\xi ')$, associated
with the orthonormal system 
of eigenvectors $e_1(x',\xi '),\dots, e_N(x',\xi ')$ in $\Bbb C^N$,
and denote$$ 
a_i(x',\xi ')=\sqrt{\lambda _i(x',\xi ')}
$$
(all equal to $\sqrt c$, when ${p'}^0=cI$). 
Denote $\Gamma _{\theta ,\pm}=\{\mu \in \Gamma _\theta \mid
\operatorname{Im}\mu \in 
\crpm\}$.
When $\mu $ runs through
$\Gamma _{\theta ,\pm}\cup\{0\}$, 
then
$(a_1^2+\mu ^2)^{\frac12}$ runs through a convex subset $V_{a_1,\pm}$ of
$\Gamma _{\theta ,\pm}$ lying to the right of a curve $C_{a_1,\pm}$
passing through $a_1$ on the real axis. 

Since $a_j\ge a_1$ for $j\ge 1$, we also have
that $(a_j^2+\mu ^2)^\frac12$ lies in
$V_{a_1,\pm}$ when $\mu \in \Gamma _{\theta ,\pm}\cup\{0\}$. We
denote $V_{a_1,+}\cup V_{a_1,-}=V_{a_1}$. It is important that although
$V_{a_1}$ is not in general convex, the $V_{a_1,\pm}$ are so.

Let $\mu \in \Gamma _{\theta ,+}\cup\{0\}$. Then for
general $v\in\Bbb C^N$ with norm 1, decomposed as
$v=c_1e_1+\dots+c_Ne_N$
with $|c_1|^2+\dots+|c_N|^2=1$,$$\aligned
\frak a^0v\cdot\bar v&=({p'}^0+\mu ^2)^\frac12 v\cdot \bar v\\
&=\sum_{i,j=1}^N(a_i^2+\mu ^2)^\frac12 c_ie_i\cdot\overline{c_je_j}=
\sum_{i=1}^N(a_i^2+\mu ^2)^\frac12 |c_i|^2\in V_{a_1,+},
\endaligned$$
since $V_{a_1,+}$ is convex. There is a similar argument
for $\mu \in \Gamma _{\theta ,-}$, showing altogether that$$
({p'}^0+\mu ^2)^\frac12 v\cdot \bar v
\in V_{a_1},\text{ when }\mu \in \Gamma _\theta \cup\{0\},\, |v|=1.\tag2.37
$$ 
Now (2.24) is obtained, if at each $(x',\xi ')$ with $|\xi '|=1$,$$
|(({p'}^0+\mu ^2)^\frac12 -\pi ^h_2b^h\pi ^h_2)v\cdot\bar v|\ge \delta|v|^2,\quad
v\in\Bbb C^N,\tag2.38
$$
for some $\delta>0$. In view of (2.37), this holds if $\pi ^h_2b^h\pi
^h_2v\cdot\bar v$ lies in a set in $\Bbb C$ with distance $\delta$
from $V_{a_1}$ when $|v|=1$.

We list some special cases where this holds; here we assume that
$\pi ^h_2$ is an 
orthogonal projection.
\flushpar {\bf (1)} Let $b^h$ be the principal symbol of a scalar first-order
differential operator with real coefficients. Then $b^h$ is purely imaginary
and$$
\pi ^h_2b^h\pi ^h_2v\cdot\bar v=b^h|\pi ^h_2v|^2\in i\Bbb R,\tag 2.39
$$
which certainly has positive distance from $V_{a_1}$.  More
generally, we can take $b^h$ such that $b^hv\cdot \bar v$ ranges in
the sectors around $i\Bbb R$ consisting of complex numbers with
argument in $\,]\theta _1,\frac\pi 2-\theta _1[\,$ or  $\,]\pi +\theta _1
,\frac{3\pi }2-\theta _1[\,$ for some $\theta _1\in \,]\theta ,\frac\pi 2]$
(allowing also pseudodifferential choices).
\flushpar {\bf (2)} Let $b^h=ib_1$, where $b_1$ is the principal symbol of a scalar
first-order 
differential operator with real coefficients. Then $b^h$ is real, and
the real number $\pi ^h_2b^h\pi ^h_2v\cdot\bar v=b^h|\pi ^h_2v|^2$
has positive distance from $V_{a_1}$ for $|v|=1$ if$$
|b^h|<a_1.\tag2.40
$$
More generally, we can take $b^h$ real selfadjoint with numerical
range in $\,[-a_1+\varepsilon ,+\infty [\,$ for some $\varepsilon \in \,]0, a_1[\,$, i.e.,$$
b^hv\cdot \bar v\ge -a_1+\varepsilon \text{ for }|v|=1.\tag2.41
$$
When $B$ is a differential operator, this in fact requires that
$|b^hv\cdot\bar 
v|\le a_1-\varepsilon $, since the symbol $b^h$ is odd in $\xi '$. This case
seems to have an interest in brane theory according to \cite{V01}.
\flushpar {\bf (3)} Let $r_1=\operatorname{dist}(V_{a_1}, 0)$. Then it
suffices that the norm of $b^h$ is $<r_1$.
\endexample

Since $R_T$ is of order $-2$ and $X$ is compact, the powers $R_T^m$
are trace-class when $m>\frac n2$. They can be studied by composition
or by differentiation, in view of the fact that $\partial _\lambda
^{m-1}(P_T-\lambda )^{-1}=(m-1)!(P_T-\lambda )^{-m}$.

\proclaim{Corollary 2.12} Under the hypotheses of Theorem {\rm 2.10},
the resolvent powers $R_T^m$ have the structure, for any $m\ge 1$:
$$
R^m_T=(Q^m)_++G^{(m)}=\tfrac1{(m-1)!}{\partial _\lambda ^{m-1}}
R_T
=\tfrac1{(m-1)!}({\partial _\lambda ^{m-1}}Q)_+ +
\tfrac1{(m-1)!}{\partial _\lambda ^{m-1}}G,
\tag2.42$$
where 
$G^{(m)}=\frac1{(m-1)!}\partial _\lambda ^{(m-1)}G$ is a singular
Green operator of class $0$  
lying in \linebreak$\operatorname{OPG}\Cal S^{0,0,-2m-1}(\Gamma _{\theta },\Cal
S_{++})$.
\endproclaim
 
\demo{Proof} One proof consists of applying the rules of calculus
(\cite{G01}, \cite{G02}) 
to the compositions $(Q_++G)\dots(Q_++G)$. Another proof is to use
the exact formula we found in Theorem 2.10, combined with the fact
that all the factors have the property that a differentiation with
respect to $\lambda $ lowers the $s$-index by 2.
\qed
\enddemo 

\proclaim{Theorem 2.13} Assumptions as in Theorem {\rm 2.10}.

{\rm (i)} Let $\varphi $ be a morphism in $E$ 
and let $m>\frac{n}2$. Then $\varphi R^m_T(\lambda )$ is trace-class and
the trace has an expansion for $|\lambda |\to\infty $ with $\arg
\lambda \in \,]\pi -2\theta , \pi +2\theta [\,$ (uniformly in closed
subsectors): $$
\Tr \bigl(\varphi R_T^m(\lambda )\bigr)
\sim
\sum_{-n\le k<0} \tilde a_{ k}(\varphi )(-\lambda ) ^{-\frac{ k}2-m}+ 
\sum_{k\ge 0}\bigl({ \tilde a'_{ k}(\varphi )}\log (-\lambda ) +{\tilde a''_{
k}(\varphi )}\bigr)(-\lambda ) ^{-\frac{ k}2-m}.\tag2.43
$$
The coefficients $\tilde a_k$ and $\tilde a'_k$ are locally
determined.

{\rm (ii)}
Let $F$ be a differential operator in $E$ of order
$m'$ and let $m>\frac{n+m'}2$. Then $FR^m_T(\lambda )$ is trace-class and
the trace has an expansion for $|\lambda |\to\infty $ with $\arg
\lambda \in \,]\pi -2\theta , \pi +2\theta [\,$ (uniformly in closed
subsectors):$$
\Tr\bigl( F R_T^m(\lambda )\bigr)
\sim
\sum_{-n\le k<0 } \tilde a_{ k}(F)(-\lambda ) ^{\frac{{m'} -k}2-m}+ 
\sum_{k\ge 0}\bigl({ \tilde a'_{ k}(F)}\log (-\lambda ) +{\tilde a''_{
k}(F)}\bigr)(-\lambda ) ^{\frac{{m'} -k}2-m},\tag2.35
$$
with locally determined coefficients $\tilde a_k$ and $\tilde a'_k$.
If
$m'$ is odd, 
$\tilde a_{-n}=0$.

Here, if $F$ is tangential (differentiates only with respect to $x'$)
on $X_c$, the log-coefficients $\tilde a'_k$ with $0\le k<m'$ vanish,
and the $\tilde 
a''_k$ with $0\le k<m'$ are locally determined.
\endproclaim

\demo{Proof} (i) is the special case of (ii) where $m'=0$, so we can
treat them at the same time. It is well-known that the ``interior''
$\psi $do term  $FQ^m_+$ produces
a series of powers $\sum_{k= -n}^\infty \tilde c_k(-\lambda )^{\frac{m'-k}
2-m}$ for $\lambda \to \infty $ in closed subsectors of $\Bbb
C\setminus \Bbb R_+$, with coefficients determined from the
successive homogeneous 
terms in the 
symbol; here the terms with  $k-n-m'$ odd vanish 
since they are produced by terms that are
odd in $\xi $. 

Now consider the singular Green term
$FG^{(m)}$.  By Corollary 2.12 and the composition rules it lies in
$\operatorname{OPG}\Cal S^{0,0,m'-2m-1}(\Gamma _\theta )$ in general,
and in $\operatorname{OPG}\Cal S^{m',0,-2m-1}(\Gamma _\theta )$ if
$F$ is tangential. With a cut-off function $\chi $ as after (2.31), we have
that $\Tr(FG^{(m)}-\chi FG^{(m)}\chi )$ is $O(|\mu | ^{-M})$
for $\mu  \to \infty $ in closed subsectors of $\Gamma _\theta $, all
$M$; this 
difference does not contribute to the trace expansion, so it 
suffices to treat $\chi FG^{(m)}\chi $. This operator can be considered as a
singular Green operator on $X^0_+=X'\times\Bbb R_+$, with the same
symbol-kernel estimates as indicated above. Now we take the
normal trace, obtaining an operator$$
\Cal S=\tr_n (\chi FG^{(m)}\chi )\in
\operatorname{OP}'S^{0,0,m'-2m}(\Gamma _\theta )\subset 
\operatorname{OP}'(S^{m'-2m,0,0}(\Gamma _\theta )\cap
S^{0,m'-2m,0}(\Gamma _\theta ))\tag 2.45 
$$
by (2.29) and (2.1). If $F$ is tangential, the last indication is
replaced by \linebreak$\operatorname{OP}'(S^{m'-2m,0,0}(\Gamma _\theta )\cap
S^{m',-2m,0}(\Gamma _\theta ))$.

Here we apply \cite{GS95, Th\. 1.2}. It is an important point in that
theorem that one gets a sum of two expansions$$
\sum _{j\ge 0}c_j\mu ^{M_0-j}+\sum_{k\ge 0}(c'_k\log\mu +c''_k)\mu
^{M_1-k} , 
\tag 2.46$$
 where $M_0$ equals the dimension (here $n-1$) plus the order of
the operator (here $m'-2m$), and $M_1$ equals the lowest $d$-index
associated with the operator. The coefficients $c_j$ and $c'_k$ are
local (each stems from a specific homogeneous term in the symbol),
whereas the $c''_k$ are global in the sense that they depend on the
full structure. So in the application to $\Cal S$, the first series
begins with a 
term $c\mu ^{n-1+m'-2m}$, and the second series begins with a term
$(c'\log\mu +c'')\mu ^{m'-2m}$ for general $F$, 
$(c'\log\mu +c'')\mu ^{-2m}$ for tangential $F$. Adding the
expansions and inserting $\mu =(-\lambda )^{\frac12}$, we find the
expansion (2.35) ff.\ after an adjustment of the indexation. 
\qed
\enddemo 

\example{Remark 2.14}
One has in particular that when all the occurring operators are strongly
polyhomogeneous, there
is only the 
first expansion in (2.46), no logarithmic or global terms.
This applies to the special case where $\Pi _1$ is a morphism in $E'$
and $B$ is a
differential operator; then$$
\Tr\bigl( F R_T^m(\lambda )\bigr)
\sim
\sum_{ k\ge -n} \tilde a_{ k}(F)(-\lambda ) ^{\frac{{m'} -k}2-m}
.\tag2.47
$$
\endexample

\subhead 3. Heat operators and power operators \endsubhead

The heat operator associated with $P_T$ is the solution operator
$u_0(x)\mapsto u(x,t)=e^{-tP_T}u_0$ for the problem$$\aligned
\partial _tu+Pu&=0\text{ on } X\times\rp,\\
\Pi _1\gamma _0u&=0\text{ on }X'\times \rp,\\
\Pi _2(\gamma _1u+B\gamma _0u)&=0\text{ on }X'\times \rp,\\
u|_{t=0}&=u_0\text{ on }X.
\endaligned\tag3.1$$
When $\theta >\frac\pi 4$ in the above constructions, the heat
operator can be defined from the resolvent powers or derivatives
(recall (2.42))
 by the formula$$
e^{-tP_T}=t^{-m}\tfrac{ i}{2\pi }\int_{{\Cal C'}}e^{-t\lambda }
\partial _\lambda ^{m}(P_T-\lambda )^{-1}
\,d\lambda;\tag3.2
$$
here $\Cal C'$
is a positively oriented curve in $\Bbb C$ going around the spectrum (like
the boundary of $\Sigma _{\delta ,R}$ in (1.9) with $\delta \in \,]0,\frac\pi
2[\,$; one can take $\delta =\pi -2\theta '$ for a $\theta '\in
\,]\frac\pi 4,\theta [\,$).

One could construct the heat operator directly instead of passing via
the
resolvent as we did above; one advantage of our approach is that we
can compose our $\lambda $-dependent operators pointwise in $\lambda
$, whereas calculations with respect to the time-variable $t$ need
convolutions. (The passage from the $\lambda $-framework to the
$t$-framework is essentially an inverse Laplace transformation; here products
are turned into convolutions, as is usual for such integral transforms.)

As shown e.g\. in \cite{GS96} (or see Sect\. 2 of \cite{G97}), the transition
formula likewise applies to 
the trace expansions, carrying the expansions in Theorem 2.13 over to
heat trace expansions with logarithms. In the resulting statement, 
we repeat our hypotheses for the convenience of the
reader:

\proclaim{Corollary 3.1} Let $P_T$ be the realization of $P$ defined
by the boundary condition {\rm (1.3)}; let
Assumptions {\rm 1.1} and {\rm 2.4} hold, and let
Assumption {\rm 2.7} hold with $\theta >\frac\pi 4$.
Then the heat operator $e^{-tP_T}$ is well-defined and its trace
has the asymptotic expansion for $t\to 0+$, for any morphism $\varphi
$ in $E $:
$$
\Tr \bigl(\varphi e^{-tP_T}\bigr)
\sim
\sum_{-n\le k<0}  a_{ k}(\varphi )t ^{\frac{ k}2}+ 
\sum_{k\ge 0}\bigl({  -a'_{ k}(\varphi )}\log t +{ a''_{
k}(\varphi )}\bigr)t^{\frac{ k}2}.\tag3.4
$$
The coefficients $ a_k$ and $ a'_k$ are locally
determined.

Moreover, when $F$ is a differential operator in $E$ of order
$m'$ there is a trace expansion for $t\to 0+$:$$
\Tr\bigl( F e^{-tP_T}\bigr)
\sim
\sum_{-n\le k<0 }  a_{ k}(F)t ^{\frac{k-{m'}}2}+ 
\sum_{k\ge 0}\bigl({ -a'_{ k}(F)}\log t +{ a''_{
k}(F)}\bigr)t ^{\frac{k-{m'}}2},\tag3.5
$$
the coefficients $ a_k$ and $ a'_k$ being locally
determined. If
$m'$ is odd, 
$ a_{-n}=0$.

Here, if $F$ is tangential 
on $X_c$, the log-coefficients $ a'_k$ with $0\le k<m'$ vanish,
and the $
a''_k$ with $0\le k<m'$ are locally determined.

The coefficients $a_k, a'_k, a''_k$ are related to the coefficients
$\tilde a_k,\tilde a'_k,\tilde a''_k$ in Theorem {\rm 2.13} by
universal nonzero proportionality factors; in particular, $$
a'_0(F)=\tilde a'_0(F)\text{ and } a''_0(F)=\tilde a''_0(F).\tag3.6
$$

\endproclaim 

Observe that as in Remark 2.14, the expansion simplifies to 
$$
\Tr\bigl( F e^{-tP_T}\bigr)
\sim
\sum_{ k\ge -n }  a_{ k}(F)t ^{\frac{k-{m'}}2},\tag3.7
$$
when $\Pi _1$ is a morphism and $B$ is a differential operator.
\medskip
The power function $FP_T^{-s}$ (defined as $0$ on
the nullspace of $P_T$) is derived from
the resolvent by the formula$$
P_T^{-s}=\tfrac i{2\pi (s-1)\dots(s-m)}\int_{\Cal C}\lambda
^{m-s}\partial _\lambda ^m(P_T-\lambda )^{-1}\,d\lambda ,\tag3.8
$$
where $\Cal C$ is a curve in $\Bbb C\setminus\crm$ around the nonzero
spectrum; here we do not need 
the extra hypothesis $\theta >\frac\pi 4$. (Further details on
transition formulas are found e.g\. in \cite{GS96} 
or \cite{G97}.)

Then one can also deduce from Theorem 2.13  that there is the
following pole structure of the
zeta function $\zeta (F,P_T,s)=\Tr(FP_T^{-s})$ for
$s\in\Bbb C$:
$$\gather
\Gamma (s)\Tr \bigl(\varphi P_T^{-s}\bigr)
\sim
\sum_{-n\le k<0}  \frac{a_{ k}(\varphi )}{s+{\frac{ k}2}}-\frac{\Tr
\varphi \Pi _0(P_T)}s + 
\sum_{k\ge 0}\bigl(\frac{ a'_{ k}(\varphi )}{(s+\frac
k2)^2} +\frac{ a''_{
k}(\varphi )}{s+\frac k2}\bigr),\tag3.9\\
\Gamma (s)\Tr\bigl( F P_T^{-s}\bigr)
\sim
\sum_{-n\le k<0 }  \frac{a_{ k}(F)}{s+\frac{k-{m'}}2}-\frac{\Tr F\Pi
_0(P_T)}s + 
\sum_{k\ge 0}\bigl(\frac{ a'_{
k}(F)}{(s+\frac{k-{m'}}2)^2} +\frac{ a''_{
k}(F)}{s+\frac{k-{m'}}2} \bigr);
\tag3.10\endgather
$$
here $a_k$, $a'_k$ and $a''_k$ are derived from $\tilde a_k$, $\tilde
a'_k$ and $\tilde a''_k$ by {\it the same} universal formulas as in
Corollary 3.1. In particular, (3.6) holds.

Consider the case where $F$ is a first-order operator $D_1$ with the
structure $D_1=$\linebreak$\psi (\partial _{x_n}+B_1)$ on $X'$ for some morphism
$\psi $. Here we get
the generalized eta-function $\Tr\bigl( D_1 P_T^{-s}\bigr)$, 
which has the pole structure
$$
\Gamma (s)\Tr\bigl( D_1 P_T^{-s}\bigr)
\sim
\sum_{-n\le k<0 }  \frac{a_{ k}(D_1)}{s+\frac{k-{1}}2}-\frac{\Tr D_1\Pi
_0(P_T)}s + 
\sum_{k\ge 0}\bigl(\frac{ a'_{
k}(D_1)}{(s+\frac{k-{1}}2)^2} +\frac{ a''_{
k}(D_1)}{s+\frac{k-{1}}2} \bigr).
\tag3.11
$$
With $s=\frac{s'+1}2$ this takes the more customary form
$$
\multline
\Tr\bigl( D_1 P_T^{-\frac{s'+1}2}\bigr)
\sim
\frac 1{\Gamma (\frac{s'+1}2)}\Bigl(\sum_{-n\le k<0 }  \frac{2a_{
k}(D_1)}{s'+k}-\frac{2\Tr D_1\Pi 
_0(P_T)}{s'+1}  \\+
\sum_{k\ge 0}\bigl(\frac{ 4a'_{
k}(D_1)}{(s'+k)^2} +\frac{ 2a''_{
k}(D_1)}{s'+k} \bigr)\Bigr).
\endmultline\tag3.12
$$
Note that since $\Gamma (\frac{s'+1}2)^{-1}$ is regular at $s'=0$,
this function in general has a double pole at $s'=0$. Moreover,
$a_{-n}(D_1)=0$ since the principal interior symbol is odd in $\xi $.

\subhead 4. The first log-term and nonlocal term\endsubhead

The hypotheses of Theorem 2.10 are assumed
throughout this section. 

It is important to investigate whether the first log-term $a'_0$ in
(3.4), resp\. $\tilde a'_0$ in (2.43) vanishes, equivalently (cf\.
(3.9)) whether 
the zeta function $\zeta (\varphi ,P_T,s)$ is regular at $s=0$.
The values of $\tilde a''_0$ and $a''_0$ are of great interest too,
and similar questions can be asked with $\varphi $ replaced by $F$,
in particular for (3.12).

In view of (3.6), we can interchange $\tilde a'_0(F)$ and $\tilde a''_0(F)$
freely with $a'_0(F)$ resp\. $a''_0(F)$.

The coefficient $a'_0$ is known to vanish in the Atiyah-Patodi-Singer
case described 
in Example 1.2, where $\Pi _1=\Pi _\ge$, $B=A_1(0)$ and $\varphi =I$, by a
slightly tricky argument (\cite{G92} proved it using \cite{APS75}; see also \cite{GS96,
pf\. of Cor\. 2.3}). 
It also vanishes in cases of local boundary 
conditions, as
mentioned in Remark 2.14.

As indicated in the proof of Theorem 2.13, the first power where
log-terms and global coefficients appear is determined from the possible
values of $d$, the second upper index in the symbol spaces $S^{m,d,0}$ that
the normal trace of the singular Green operator belongs to. 
Consider first the case where $F=\varphi $, a morphism in $E$ (taken
to be
constant in $x_n$ on $X_c$). Then
Theorem 2.13 tells us that the index $d$ can be taken equal to $-2m$.
We shall in 
the following strive to isolate the part of $G^{(m)}$ that
contributes nontrivially to the first log-term and global term, in
a way that gives information on the value.

\proclaim{Theorem 4.1} Let $\varphi $ be a morphism in $E$,
independent of $x_n$ on $X_c$ (its restriction to $X'$ likewise
denoted $\varphi $).  Consider $\chi \varphi G^{(m)}\chi $ (cf\. Theorem
{\rm 2.13}) as an operator on
$X^0_+=X'\times \Bbb R_+$. With the notation introduced around {\rm
(2.32)}, we have that
$$\aligned
\chi \varphi G^{(m)}\chi &=\tfrac{\partial _\lambda
^{m-1}}{(m-1)!}\varphi  \Pi _2K_{\frak A}\frak
A^{-1}T_{\frak A}+\Cal G_1+\Cal G_2,\\
\tr_n(\chi \varphi G^{(m)}\chi )&=\tfrac{\partial _\lambda
^{m-1}}{(m-1)!}\tfrac12\varphi \Pi _2\frak 
A^{-2}+\Cal S_1+\Cal S_2,
\endaligned\tag 4.1$$
where $\Cal G_1$ is a strongly polyhomogeneous singular Green
operator of degree $-2m-1$, $\Cal G_2\in \operatorname{OPG}\Cal
S^{1,0,-2-2m}(\Gamma _\theta ,\Cal S_{++})$, $\Cal S_1$ is a strongly
polyhomogeneous $\psi $do in $E'$ of 
degree $-2m$, and  $\Cal S_2\in \operatorname{OP}'S^{1,0,-1-2m}(\Gamma
_\theta )$. 
\endproclaim 

\demo{Proof} By Theorem 2.10 and (2.42), we have that $$
\chi \varphi G^{(m)}\chi =-\tfrac{\partial _\lambda
^{m-1}}{(m-1)!}\chi \varphi K_{\operatorname{D}}\gamma
_0Q_+\chi +\tfrac{\partial _\lambda ^{m-1}}{(m-1)!}\chi
\varphi K_{\operatorname{D}} [S_0\gamma 
_0+S_1\gamma _1]Q_+\chi .
$$
The first term is a $\Cal G_1$. For the second term, we
first note that, since $(1-\chi )\partial _\lambda ^rK_{\frak A}$ and
$\partial _\lambda ^rT_{\frak A}(1-\chi
)$ are smoothing and $O(|\mu |^{-M})$ in closed subsectors of $\Gamma
$ for all $M$ (cf\. e.g\. \cite{G01,
Lemma 7.1}),
$$\aligned
\partial _\lambda ^r\chi K_{\operatorname{D}}-\partial _\lambda ^rK_{\frak
A}&\in \operatorname{OPK}\Cal S^{0,0,-2-2r}(\Gamma ,\Cal S_+),\\
\partial _\lambda ^r\gamma _0Q_+\chi-\partial _\lambda ^r\tfrac12\frak
A^{-1}T_{\frak A }&\in  
\operatorname{OPT}\Cal S^{0,0,-3-2r}(\Gamma ,\Cal S_+),\\  
\partial _\lambda ^r\gamma _1Q_+\chi-\partial _\lambda ^r\tfrac12
T_{\frak A }&\in  
\operatorname{OPT}\Cal S^{0,0,-2-2r}(\Gamma ,\Cal S_+),  
\endaligned\tag 4.2$$ 
in view of the information around (2.32)--(2.34).
Thus, using the formulas for $S_0$, $S_1$ and $S$ and replacing
$A_{\operatorname{DN}}$ by $-\frak A$, 
$$
\aligned
\chi \varphi G^{(m)}\chi &=\tfrac{\partial _\lambda
^{m-1}}{(m-1)!}\varphi K_{\frak A}S^{-1}\Pi _2[(-\frak
A)\tfrac12\frak A^{-1}-\tfrac12]T_{\frak A}+\Cal G_1+\Cal G'_2,\\ 
&=-\tfrac{\partial _\lambda
^{m-1}}{(m-1)!}\varphi K_{\frak A}(-\frak A+\Pi _2B\Pi _2)^{-1}\Pi _2T_{\frak A}+\Cal G_1+\Cal G''_2,
\endaligned\tag4.3$$
where $\Cal G'_2$ and $\Cal G''_2$ are in $\operatorname{OPG}\Cal
S^{0,0,-2-2m}(\Gamma _\theta ,\Cal S_{++})$.
We have furthermore:$$
(-\frak A+\Pi _2B\Pi _2)^{-1}\Pi _2
=-\frak A^{-1}\Pi _2+(-\frak A+\Pi _2B\Pi _2)^{-1}\Pi
_2B\Pi _2\frak A^{-1}\Pi _2
.\tag4.4$$
Here $\partial _\lambda ^{m-1}$ of 
$\frak A^{-1}\Pi _2$ is in
$\operatorname{OP}'S^{0,0,1-2m}(\Gamma )$. Since  
$\partial _\lambda ^r(-\frak A+\Pi _2B\Pi _2)^{-1}$ is in \linebreak$
\operatorname{OP}'S^{0,0,-1-2r}(\Gamma _\theta )$, $\Pi
_2B\Pi _2\in \operatorname{OP}'S^{1,0,0}(\Gamma )$ and $\partial
_\lambda ^r\frak A^{-1}\Pi _2\in \operatorname{OP}'S^{0,0,-1-2r}(\Gamma )$,
 $\partial _\lambda ^{m-1}$ of the last term is
in $\operatorname{OP}'S^{1,0,-2m}(\Gamma _\theta )$. Thus,
if we use Proposition 2.5 to replace $\partial _\lambda ^r\frak
A^{-1}\Pi_2$ by 
$\partial _\lambda ^r\Pi_2\frak A^{-1}$, we can write 
$$
\chi \varphi G^{(m)}\chi =\tfrac{\partial _\lambda
^{m-1}}{(m-1)!}\varphi  K_{\frak A}\Pi _2\frak
A^{-1}T_{\frak A}+\Cal G_1+\Cal G'''_2,\tag 4.5
$$
with 
$\Cal G'''_2$ like $\Cal G_2$ in
the theorem. Note that $B$ has disappeared from the main term!

It will be convenient to do one more commutation, placing $\Pi_2$ in
front of the Poisson operator $K_{\frak
A}$, which likewise gives an error like $\Cal G_2$. For this we use
that $e^{-x_n\frak A}=$\linebreak$\operatorname{OPK}_n((\frak A+i\xi
_n)^{-1})$ 
in the complex formulation.
Here 
$[(\frak A+i\xi _n)^{-1},\Pi _2]=$\linebreak$(\frak A+i\xi _n)^{-1}[\Pi _2,\frak
A](\frak A+i\xi _n)^{-1}$, which by application of
$\operatorname{OPK}_n$ gives a Poisson operator in
$\operatorname{OPK}\Cal S^{0,0,-2}(\Gamma _\theta ,\Cal S_+)$ (with
$(m-1)$'st derivative in $\operatorname{OPK}\Cal S^{0,0,-2m}(\Gamma
_\theta ,\Cal S_+)$), so 
that we can write
$$
\chi  \varphi G^{(m)}\chi =\tfrac{\partial _\lambda
^{m-1}}{(m-1)!}\varphi \Pi _2K_{\frak A}\frak
A^{-1}T_{\frak A}+\Cal G_1+\Cal G^{(4)}_2,
$$
with $\Cal G_2^{(4)}$ like $\Cal G_2$.
This shows the first formula in (4.1). 

By functional calculus as in (2.34), we have
that $\tr_n (K_{\frak A}(\frak A)^{-1}T_{\frak A})=\frac12 \frak
A^{-2}$. Then also $\tr_n (\varphi \Pi _2K_{\frak A}(\frak
A)^{-1}T_{\frak A})= \frac12\varphi \Pi _2 \frak 
A^{-2}$, showing the second formula in (4.1); the symbol properties
follow from (2.29). (Placing the $\psi$do
$\varphi \Pi _2$ in front of the singular Green operator $K_{\frak 
A}(\frak A)^{-1}T_{\frak A}$, and taking it outside $\tr_n$, is
justified by the point of view where $\operatorname{OPG}=\operatorname{OPG}_n\operatorname{OP}'$.) 
\qed 
\enddemo 

The last commutation in the proof is not needed when $\varphi =I$,
for then we can instead use circular permutation in the treatment of
$\tfrac{\partial _\lambda
^{m-1}}{(m-1)!} K_{\frak A}\Pi _2\frak
A^{-1}T_{\frak A}$:$$
\Tr_{X^0_+}(\tfrac{\partial _\lambda
^{m-1}}{(m-1)!}K_{\frak A}\Pi _2\frak
A^{-1}T_{\frak A})=\Tr_{X'}(\tfrac{\partial _\lambda
^{m-1}}{(m-1)!}\Pi _2\frak
A^{-1}T_{\frak A} K_{\frak A})=
\tfrac12\Tr_{X'}(\tfrac{\partial _\lambda
^{m-1}}{(m-1)!}\Pi _2\frak
A^{-2}),\tag 4.6
$$
since $T_{\frak A} K_{\frak A}=(2\frak A)^{-1}$.

Observe that
$$
\varphi \Pi _2\tfrac{\partial _\lambda ^{m-1}}{(m-1)!}   \frak A^{-2}=
\varphi \Pi _2\tfrac{\partial _\lambda ^{m-1}}{(m-1)!} (P'-\lambda
)^{-1};\tag4.7 
$$
it lies in $\operatorname{OP}'S^{0,0,-2m}(\Gamma )$.
Let us list the three trace expansions connected with this operator
family (which follow from \cite{GS95, Th\. 2.7}):$$\align
\Tr_{X'}(\tfrac12\varphi \Pi _2\tfrac{\partial _\lambda
^{m-1}}{(m-1)!}&(P'-\lambda )^{-1})\sim 
\sum_{ k\ge n-1} \tilde c_{ k}(\varphi )(-\lambda ) ^{-\frac{
k}2-m}\\
&\qquad+ 
\sum_{l\ge 0}\bigl({ \tilde c'_{ 2l}(\varphi )}\log (-\lambda )
+{\tilde c''_{2 l}(\varphi )}\bigr)(-\lambda ) ^{-l-m},\\
\Tr_{X'}(\tfrac12\varphi \Pi _2e^{-tP'})&\sim 
\sum_{ k\ge 1-n}  c_{ k}(\varphi )t ^{\frac{
k}2}+ 
\sum_{l\ge 0}\bigl({- c'_{2 l}(\varphi )}\log t
+{ c''_{2l}(\varphi )}\bigr)t ^{l},\tag4.8\\
\Gamma (s)\Tr_{X'}(\tfrac12\varphi \Pi _2(P')^{-s})&\sim 
\sum_{ k\ge 1-n} \frac{c_{ k}(\varphi )}{s+\frac k2}-\frac{\Tr[\frac12\varphi
\Pi _2\Pi _0(P')]}s
\sum_{l\ge 0}\bigl(\frac{  c'_{2l}(\varphi )}{(s+l)^2}
+\frac{ c''_{2l}(\varphi )}{s+l}\bigr);
\endalign$$
again with $c'_0(\varphi )=\tilde c'_0(\varphi )$, $c''_0(\varphi )=\tilde c''_0(\varphi )$.

\proclaim{Theorem 4.2} Let $\varphi $ be as in Theorem {\rm 4.1}. 
In a comparison of {\rm (2.43)} with the trace expansion of
$\tfrac12\varphi \Pi _2\frac{\partial _\lambda 
^{m-1}}{(m-1)!}(P'-\lambda )^{-1}$ in {\rm (4.8)},
we have that$$\aligned
\tilde a'_0(\varphi )&=\tilde c'_0(\varphi ),\\
\tilde a''_0(\varphi )&=\tilde c''_0(\varphi )+\text{ local contributions,}
\endaligned\tag 4.9$$
where the local contributions are defined from $P$ (at $X'$), $\varphi $ and
the strictly homogeneous terms in the symbol of $\Pi _2$ of orders
$\{0,-1,\dots,1-n\}$.

In particular, $\tilde a'_0(\varphi )$
equals the residue $c'_0(\varphi )$ at $0$ of the zeta
function $\zeta (\tfrac12\varphi \Pi _2,P',s)=\Tr_{X'}(\tfrac12\varphi
\Pi _2(P')^{-s})$, also identifiable with the non-commutative residue
of $\tfrac12\varphi \Pi _2$ times $\tfrac12$:$$ 
\tilde a'_0(\varphi )=a'_0(\varphi )=c'_0(\varphi
)=\operatorname{Res}_{s=0}\Tr(\tfrac12\varphi \Pi _2(P')^{-s})
=\tfrac14\operatorname{res}(\varphi \Pi _2).\tag4.10
$$
Moreover, when $\theta >\frac\pi 4$, $a'_0(\varphi )$ in {\rm (3.4)}
equals the 
coefficient $c'_0(\varphi )$ of $\log t$ in the trace expansion of
$\tfrac12\varphi \Pi 
_2e^{-tP'}$
in {\rm(4.8)}.
\endproclaim

\demo{Proof} As noted in Theorem 2.13, the part $\varphi (Q^m)_+$ of
$\varphi R_T^m$
contributes only a local power expansion, and $\varphi G^{(m)}$
contributes the same expansion as $\chi \varphi G^{(m)}\chi $, so it
suffices to study 
$\Tr_{X'}\tr_n 
(\chi \varphi G^{(m)}\chi )$, where we can use (4.1). Here the
strongly polyhomogeneous part $\Cal S_1$ gives an expansion 
with purely local power terms (as in (2.47) with $k\ge 1-n$ and
$m'=0$), and the term $\Cal S_2$ with symbol in $S^{1,0,-1-2m}(\Gamma
_\theta )\subset
S^{-2m,0,0}(\Gamma _\theta )\cap S^{1,-1-2m,0}(\Gamma _\theta )$ gives
an expansion as in (2.46) with  
$M_0=n-2m$ and $M_1=-1-2m$.
So the non-local and log-contributions from $\Cal S_1+\Cal S_2$ begin
with a term $(c'\log(-\lambda )+c'')(-\lambda )^{-\frac12-m}$,
implying (4.9).
 
The identification of $c'_0(\varphi )$ with $\tfrac12$ times the
noncommutative 
residue $\operatorname{res}(\tfrac12\varphi \Pi _2)$ refers to the work of
Wodzicski \cite{W84}, \cite{W84$'$}, where it was shown that for a classical
integer-order $\psi $do $C$ in $E'$, the following formula:$$
\operatorname{res}(C)=\operatorname{ord}Q\,
\operatorname{Res}_{s=0}\Tr (CQ^{-s} )
=\int_{X'}\int_{|\xi '|=1}(\tr_{E'}c(x',\xi
'))_{1-n}\,d\sigma (\xi ')\tag4.11
$$ can be given a sense (the trace $\Tr$ is defined as a meromorphic
extension from large $\operatorname{Re}s$ to $s\in \Bbb C$, $Q$ is an
auxiliary invertible elliptic $\psi $do of 
positive order, $\tr_{E'}$ is the fiber trace, and subscript $1-n$ indicates
the term of degree $1-n=-\operatorname{dim}X'$ in the
symbol); the functional res vanishes on commutators.

If $P'$ is invertible, the last equality sign in (4.10) follows
directly by taking $Q=P'$. If $P'$ is not invertible but has a
(necessarily finite 
dimensional) nullspace $V_0(P')$, we use that $(P')^{-s}$ is defined
to be 0 on that nullspace, and that a replacement of $P'$ by the
invertible operator $$
Q=P'+\Pi _0(P')\tag4.12
$$  
in (4.8) leaves the coefficient $c_0'(\varphi )$ invariant (whereas
$c''_0(\varphi )$ is changed).  
\qed
\enddemo

These formulas have consequences for the original
Atiyah-Patodi-Singer problem, that we take up in Section 5. 
Let us
also observe:

\proclaim{Theorem 4.3} Let $D_1$ be a first-order differential
operator on $X$, of the form 
$D_1=\psi (\partial _{x_n}+B_1)$ on $X_c$, where $B_1$ is
tangential and $\psi $ is a morphism independent of $x_n$. 
Then
$$\aligned
\chi D_1 G^{(m)}\chi &= -\psi \Pi _2\tfrac{\partial _\lambda
^{m-1}}{(m-1)!}K_{\frak A}T_{\frak A}+\Cal G'_1+\Cal G'_2,\\
\tr_n (\chi D_1 G^{(m)}\chi )&= -\tfrac12\psi \Pi _2\tfrac{\partial _\lambda
^{m-1}}{(m-1)!}\frak A^{-1}+\Cal S'_1+\Cal S'_2,
\endaligned\tag 4.13$$
where $\Cal G'_1$ is a strongly polyhomogeneous singular Green
operator of degree $-2m$, $\Cal G'_2\in \operatorname{OPG}\Cal
S^{1,0,-1-2m}(\Gamma _\theta ,\Cal S_{++})$, $\Cal S'_1$ is a strongly
polyhomogeneous $\psi $do in $E'$ of 
degree $-2m+1$ and  $\Cal S_2\in \operatorname{OP}'S^{1,0,-2m}(\Gamma
_\theta )$. 
Here $-\tfrac12\psi \Pi _2\tfrac{\partial _\lambda
^{m-1}}{(m-1)!}\frak A^{-1}=-\tfrac12\psi \Pi _2\tfrac{\partial _\lambda
^{m-1}}{(m-1)!}(P'-\lambda
)^{-\frac12}$ has a trace expansion:  
$$\aligned
-\Tr_{X'}(\tfrac12\psi \Pi _2\tfrac{\partial _\lambda
^{m-1}}{(m-1)!}(P'-\lambda )^{-\frac12})&\sim 
\sum_{ k\ge n-1} \tilde d_{ k}(\psi )(-\lambda ) ^{\frac{1-
k}2-m}\\
&\quad+ 
\sum_{k\ge 0}\bigl({ \tilde d'_{ k}(\psi )}\log (-\lambda )
+{\tilde d''_{k}(\psi )}\bigr)(-\lambda ) ^{\frac{1-k}2-m},
\endaligned\tag 4.14$$
where the $\tilde d'_k$ and $\tilde d''_k$ are zero for $k$ odd.
Hence in a comparison of {\rm (2.35)} for $F=D_1$ with {\rm (4.13)},
we have that$$\aligned
\tilde a'_0(D_1)&=\tilde d'_0(\psi ),\\
\tilde a''_0(D_1)&=\tilde d''_0(\psi )+\text{ local contributions,}
\endaligned\tag 4.15$$
where the local contributions are defined from $P$ and $D_1$ (at $X'$), $\psi $ and
the strictly homogeneous terms in the symbol of $\Pi _2$ of orders
$\{0,-1,\dots,1-n\}$.

(The same formulas hold for $a'_0(D_1)$ resp\. $a''_0(D_1)$.)

The coefficient $\tilde d'_0(\psi )$ is proportional to the noncommutative
residue of $\psi  \Pi _2$:
$$
\tilde d'_0(\psi )=\alpha  \operatorname{res}(\psi \Pi _2),\tag4.16
$$
for some nonzero constant $\alpha $ depending only on $m$ and the dimension.
\endproclaim 

\demo{Proof} (4.13) follows from (4.1), when we use that the
tangential operator $B_1$ lifts the first upper index in our symbol
classes by 1, $\partial _{x_n}e^{-x_n\frak A}=-\frak Ae^{-x_n\frak
A}$, and $[D_1,\chi ]$ is supported away from $x_n=0$.
Now (4.14) holds by application of \cite{GS95, Th\. 2.1} (with $\mu
^2=-\lambda $ as usual); the more precise information that the
log-terms and nonlocal terms only occur for even $k$ --- which is not
needed for our main purposes here --- follows from
\cite{GH02} (or from an analysis as in \cite{G02, Th\. 5.2}). Since $\Cal S'_1$
produces only local
power terms, and the logarithmic and global terms for $\Cal S'_2$
begin with the power $-m$, the result
(4.15) on the zero'th 
coefficients follows.

For (4.16) one applies an analysis as in \cite{G02, proof of Th\. 5.2} to
(4.14). We have that $$
\tfrac{\partial _\lambda
^{m-1}}{(m-1)!}(P'-\lambda )^{-\frac12}=c_m(P'-\lambda )^{-\frac12
-m}=(-\lambda )^{-\frac12-m}(\varrho P'+1)^{-\frac12-m},\quad \varrho
=-\lambda ^{-1}.\tag 4.17
$$
Now$$
\psi \Pi _2(\varrho P'+1)^{-\frac12-m}\sim \sum_{j \ge
0}\tbinom{-\frac 12-m}{j }\varrho ^{j }\psi \Pi _2(P')^{j
}\tag4.18
$$
similarly to \cite{G02, (5.6)}, and the subsequent analysis there gives that the
first log-term comes from the term with $j =0$ and equals 
$$
c \int_{X'}\int_{|\xi '|=1}(\tr [\psi (x')\Pi _2(x',\xi
')])_{1-n}\,d\sigma (\xi ') =c'\operatorname{res}(\psi \Pi _2) ,
\tag4.19$$
with nonzero universal constants $c$ and $c'$.
\qed
\enddemo 

The result has an
interest for the analysis of the eta function associated with the
APS-problem, which we take up in Section 5.

\example{Remark 4.4}
Let us comment on what is meant by ``local
contributions'' in (4.9) and (4.15). It is taken in a rather strict
sense, based on the derivation
of asymptotic formulas in \cite{GS95, Th\. 2.1} (also recalled in \cite{G02,
Th\. 2.10}). 
The local contributions at a certain index $k$ come from the homogeneous
symbol terms with the degree that matches the index. 

Since $\varphi R^m_T=\varphi Q^m_++\varphi G^{(m)}$ is of order $-2m$, local
contributions at $k=0$ (in (4.9)) come from terms of
homogeneity degree $-2m-n$ in $(\xi ,\mu )$ in the symbol of $\varphi
Q^m$ and terms of homogeneity degree $-2m-n+1$ in $(\xi ',\mu )$ in the
symbol of 
$ \tr_n \varphi G^{(m)}$, in local trivializations. 

Since $D_1 R^m_T=D_1 Q^m_++D_1 G^{(m)}$ is of order $-2m+1$, local
contributions at $k=0$ (in (4.15)) come from terms of
homogeneity degree $-2m+1-n$ in $(\xi ,\mu )$ in the symbol of $D_1
Q^m$ and terms of homogeneity degree $-2m-n+2$ in $(\xi ',\mu )$ in the
symbol of 
$ \tr_n D_1G^{(m)}$, in local trivializations. 

These contributions can be
traced back to the symbols of $\varphi $ and $P$, resp\. $\psi $, $P$
and $D_1$ at $x_n=0$, and the 
$x_n$-derivatives up to order $n$ at $x_n=0$ of the symbol of $P$
(resp\. $P$ and $D_1$), together
with the first $n$ homogeneous terms (down to order $1-n$) in the symbol of
$\Pi _2$ (or $\Pi _1$).

\endexample

In many cases one can show that the first
log-coefficient vanishes and get some information on the nonlocal content of
the coefficient behind it.

We first observe:

\proclaim{Theorem 4.5} We always have that $a'_0(I)=0$. Moreover,
$a'_0(\varphi )=0$ if $\varphi \Pi _2$ is a projection, and
$a'_0(D_1)=0$ if $\psi \Pi _2$ is a projection.
\endproclaim 

\demo{Proof}This follows from the fact, observed already by Wodzicki
in \cite{W84}, \cite{W84$'$}, that when $\Pi $ is a $\psi $do projection, then
$\operatorname{res}(\Pi )=0$. The three statements in the theorem
follow by applying this to the formulas shown above:$$\aligned
a'_0(I)&=c'_0(I)=\tfrac14\operatorname{res}(\Pi _2),\\
a'_0(\varphi )&=c'_0(\varphi )=\tfrac14\operatorname{res}(\varphi \Pi _2),\\
a'_0(D_1)&=d'_0(\psi )=\alpha \operatorname{res}(\psi \Pi _2);
\endaligned\tag 4.20$$
they give 0 when the operator to which res is applied is a projection.
\qed\enddemo 

Other systematic results can be obtained when 
$\Pi _1$ is a spectral
projection or a suitable perturbation of such an operator.

\proclaim{Definition 4.6}
A pseudodifferential orthogonal projection $\Pi $ in $L_2(E')$ will be said to be a
spectral 
projection associated with
$C$, when $C$ is a selfadjoint first-order elliptic classical $\psi $do
in $E'$ and $\Pi $ satisfies
$$\Pi =\Pi _>(C)+\Pi _{V'_0},
\tag4.21$$ 
where $\Pi _>(C)$ is the orthogonal projection onto the positive
eigenpace of $C$ and $\Pi _{V'_0}$ is the orthogonal projection onto
a subspace 
$V'_0$ of $V_0(C)$ (the nullspace of $C$). We denote $V_0(C)\ominus
V'_0=V''_0$. 
\endproclaim

It is shown below in Proposition 4.8 that for any orthogonal $\psi $do
projection $\Pi $ there exists an invertible $C$ such that $\Pi =\Pi
_>(C)$.

A particular result is the following:

\proclaim{Theorem 4.7}
Let 
$$\Pi _1=\Pi _>(C)+\Cal S,
\tag4.22$$
where $C$ is a selfadjoint elliptic {\bf differential} operator of order
$1$ and $\Cal S$ is a $\psi $do of order $\le -n$.
 If $n$ is odd, then $\operatorname{res}(\varphi \Pi _2)$ and
$\operatorname{res}(\psi \Pi _2)$ are zero and hence $\tilde
a'_0(\varphi )$ and 
$\tilde a'_0(D_1)$ in {\rm (4.9), (4.15)} are zero.
\endproclaim 

\demo{Proof} We have that $\varphi \Pi _2=\varphi \Pi _\le(C)-\varphi
\Cal S$, where$$ 
\Pi _\le(C)=\tfrac {|C|-C}{2|C'|}+\Pi _0(C)=\tfrac12(I-\tfrac
C{|C'|}+\Pi _0(C)),\text{ with }C'=C+\Pi _0(C). \tag4.23
$$
Here $\frac12\varphi $, $\varphi \Pi _0$ and $\varphi \Cal S$ have
residue 0, since they have no 
$1-n$-degree term in the symbol. The
symbol of $\frac C{|C'|}$ has even-odd parity (the terms of even
degree of homogeneity order are odd in $\xi '$ and vice versa; more
on such symbols e.g\. in \cite{G02, Sect\. 5}), and
so does the symbol composed with $\varphi $. Thus, when
the interior dimension $n$ is odd, the term of order $1-n$ in the
symbol is odd in $\xi '$, so the integration with respect to $\xi '$ in (4.11)
gives zero.\qed
\enddemo

Before considering the more delicate results on the term $a''_0$, we include
some words about $\psi $do projections.

\proclaim{Proposition 4.8} 
\flushpar{\rm (i)} When $\Pi $ is a $\psi $do projection in $L_2(E')$,
then
$$
\Pi _{\operatorname{ort}}=\Pi \Pi ^*[\Pi \Pi ^*+(I-\Pi ^*)(I-\Pi
)]^{-1}\tag4.24 
$$ 
is an orthogonal $\psi $do projection with the same range. Moreover,
$$
R=\Pi +(I-\Pi _{\operatorname{ort}})(I-\Pi )\tag4.25
$$
is an invertible elliptic zero-order $\psi $do (with inverse $\Pi
_{\operatorname{ort}}+(I-\Pi )(I-\Pi _{\operatorname{ort}})$) such
that $$
\Pi _{\operatorname{ort}}=R\Pi R^{-1}.\tag4.26
$$

\flushpar{\rm (ii)} 
Let $\Pi $ be an orthogonal $\psi $do projection
in $L_2(E')$. There exists a selfadjoint invertible elliptic $\psi
$do $C$ of order $1$ in $E'$ such that $\Pi =\Pi _>(C)$.
\endproclaim 

\demo{Proof} 
(i). The formula (4.24) is known from Birman and Solomyak \cite{BS82},
details of verification can also be found in Booss-Bavnbek and
Wojciechowski \cite{BW93, Lemma 12.8}. The statements on $R$ are easily
checked.

(ii). If the principal symbol $\pi ^0$ equals the identity or 0, we
are in a trivial case, so let us assume that $\pi ^0\ne I$ and $0$;
then $\Pi $ and $\Pi ^\perp$ both have 
infinite dimensional range. Let $C_1$ be a selfadjoint positive first-order elliptic
$\psi $do with scalar principal symbol $c_1^0(x',\xi ')$ (e.g\. $=
|\xi '|$). 
Let $$
C'=\Pi C_1\Pi -\Pi ^\perp C_1\Pi ^\perp.
$$
$C'$ is a $\psi $do of order 1 with principal symbol $$
{c'}^0=\pi ^0c_1^0\pi ^0-(I-\pi ^0)c_1^0(I-\pi ^0)=(2\pi ^0-I)c_1^0,
$$
which is invertible since $2\pi ^0-I$ and $c_1^0$ are so. Clearly,
$C'$ is selfadjoint, and $\Pi C_1\Pi \ge 0$, $\Pi ^\perp C_1\Pi
^\perp\ge 0$ in view of the 
positivity of $C_1$. Moreover, $\Pi $ commutes with $C'$.

 Since $C'$ is selfadjoint elliptic, it has a spectral
decomposition in smooth orthogonal finite dimensional
ei\-gen\-spa\-ces $V _k$ with mutually distinct
eigenvalues $\lambda _k$, $k\in\Bbb Z$, such that $\lambda _k<0$ for
$k<0$,
$\lambda _k>0$ for $k>0$, $\lambda _0=0$ ($V_0$ may be 0). The positive resp\.
negative eigenspace 
of $C'$ is $V_>=\oplus_{k>0}V_k$ resp\. $V_<=\oplus_{k<0}V_k$. 
The
eigenspaces are invariant under $\Pi $:
For  $u_k\in V_k$,$$
C'\Pi u_k=\Pi C'u_k=\Pi \lambda _k u_k=\lambda _k\Pi u_k;
$$ hence $\Pi u_k\in V_k$. 
Now if  $u_k\in V_k\setminus\{0\}$ with $k<0$, then $\Pi u_k$ must be
zero, for otherwise $$
(C'\Pi u_k,\Pi u_k)=\lambda _k\|\Pi u_k\|^2<0,
$$ in contrast with
$$
(C'\Pi u_k,\Pi u_k)=(\Pi C_1\Pi \Pi u_k,\Pi u_k)=(C_1\Pi u_k,\Pi u_k)\ge
c\|\Pi u_k\|^2>0;
$$ here $c>0$ is the lower bound of $C_1$. Thus $V_<\subset R(\Pi
^\perp)$. Similarly, $V_>\subset R(\Pi )$.

Finally, let $\Pi _{0+}$ be the orthogonal projection onto $\Pi V_0$,
let $\Pi _{0-}$ be the orthogonal projection onto $\Pi ^\perp V_0$
(note that $V_0=\Pi V_0\oplus \Pi ^\perp V_0$ since $\Pi V_0\subset V_0$),
and set
$$
C=C'+\Pi _{0+}-\Pi _{0-};
$$it is injective.
Then $V_>(C)=V_>(C')\oplus \Pi V_0\subset R(\Pi )$ and
$V_<(C)=V_<(C')\oplus \Pi ^\perp V_0\subset R(\Pi ^\perp)$, so since
they are complementing subspaces, they equal $R(\Pi )$ resp\. $R(\Pi
^\perp)$, and $C$ is as asserted. \qed
\enddemo

Parts of the above proof details are
given in Br\"uning and Lesch \cite{BL99, Lemma 2.6}. They can be used to
show the fact that $\operatorname{res}(\Pi )=0$ for any $\psi $do
projection, that we used above in Theorem 4.5. In fact, with the
notation of the proposition, we have since res vanishes on commutators,
$$\multline
\operatorname{res}(\Pi )=\operatorname{res}(R^{-1}\Pi
_{\operatorname{ort}}R)
=\operatorname{res}(\Pi
_{\operatorname{ort}})
=\operatorname{res}(\Pi
_>(C))=\tfrac12\operatorname{res}(I+C|C|^{-1}) 
\\=\tfrac12\operatorname{res}(C|C|^{-1}) 
=\tfrac12\operatorname{Res}_{s=0}\Tr(C|C|^{-s-1})
=\tfrac12\operatorname{Res}_{s=0}\eta
(C,s)=0,\endmultline\tag4.27
$$ 
where the last equality follows from the vanishing of the eta residue
of $C$ shown by Atiyah, Patodi 
and Singer \cite{APS76} (odd dimensions) and Gilkey \cite{Gi81}. (The
relation between the vanishing of the noncommutative residue on
projections, and the 
vanishing of eta residues, enters also in \cite{W84}.)

\proclaim{Theorem 4.9} Let $\Pi _1$ be an orthogonal
pseudodifferential projection, and let $C$ be a first-order
selfadjoint elliptic $\psi $do such that $\Pi _1=\Pi
_>(C)+\Pi _{V'_0}$ as in Definition {\rm 4.6}. Then 
$$
\Pi _2(P' -\lambda )^{-m}-\Pi _2(C^2-\lambda
)^{-m}\in \operatorname{OP}'S^{2,0,-2m-2}(\Gamma ). \tag4.28
$$
The  power function $\tfrac12\Pi _2(C^2)^{-s}$ corresponding to $\tfrac12\Pi
_2(C^2-\lambda )^{-m}$ (cf\. {\rm (3.8)}) satisfies $$\aligned
\tfrac12\Pi _2(C^2)^{-s}&=\tfrac14((C^2)^{-s}-
C|C|^{-2s-1}),\\
\Tr[\tfrac12\Pi _2(C^2)^{-s}]&=\tfrac14\zeta (C^2,s)-
\tfrac14\eta (C,2s).
\endaligned\tag4.29
$$In particular, in {\rm (4.8)},
$$
 c''_0(I)=-\tfrac14 (\eta (C,0)+\operatorname{dim}V'_0-\operatorname{dim}V''_0
)+ \text{ local contributions}.\tag4.30
$$

It follows that when $\Pi _1$ is the projection entering in the
construction in Theorems {\rm 2.10} and {\rm 2.13}, then
$$
\tilde a''_0(I)=-\tfrac14 (\eta (C,0)+\operatorname{dim}V'_0
-\operatorname{dim}V''_0)+\text{ local contributions.}
\tag 4.31$$
\endproclaim

\demo{Proof} 
It follows from \cite{GS95} (adapted to the present notation)
that $(C^2-\lambda )^{-m}$
\linebreak
$\in \operatorname{OP}'S^{0,0,-2m}(\Gamma )$,
and that
$$
\multline
\Pi _2(P' -\lambda )^{-m}-\Pi _2(C^2-\lambda
)^{-m}\\=\Pi _2
\tfrac{\partial _\lambda
^{m-1}}{(m-1)!}[(P' -\lambda )^{-1}(C^2-P')(C^2-\lambda
)^{-1}] 
\in 
\operatorname{OP}'S^{2,0,-2m-2}(\Gamma ),
\endmultline\tag4.32
$$
so this difference contributes no log-terms or nonlocal terms at the
powers $(-\lambda )^{-m}$ and $(-\lambda )^{-m-\frac12}$.
(This reflects the fact that in the residue construction, one
can replace the auxiliary operator $P'$ by $C^2$.) Now in (4.29), the first line follows from (4.23) when we recall that
the powers are defined to be zero on $V_0(C)$; the second line
follows by taking the trace. Then since $\Tr[\frac12\Pi _2\Pi
_0(C)]=\frac12\operatorname{dim}V''_0$ and $\zeta
(C^2,0)=-\operatorname{dim}V_0(C) +$ a local coefficient, we have at $s=0$:$$
 \aligned c''_0(I)-\tfrac12\operatorname{dim}V''_0&=\Tr[\tfrac12\Pi
_2(C^2)^{-s}]_{s=0}=
\tfrac14\zeta (C^2,0)-
\tfrac14\eta (C,0))\\
&=
-\tfrac14 \operatorname{dim}V_0(C)-\tfrac14\eta (C,0)+\text{ a local coefficient},
\endaligned$$
which implies (4.30) since $-\frac14 \operatorname{dim
}V_0(C)+\frac12\operatorname{dim}V''_0=-\frac14(\operatorname{dim}V'_0
-\operatorname{dim}V''_0)$.
Finally, (4.31) follows in view of (4.9).\qed
\enddemo

We can give a name to the ``nonlocal part of $a''_0(I)$'' appearing
in this way:

\proclaim{Definition 4.10} In the situation of Definition {\rm 4.6}, we
define the associated eta-invariant $\eta _{C,V'_0}$ by:
$$
\eta _{C,V'_0}=\eta
(C,0)+\operatorname{dim}V'_0-\operatorname{dim}V''_0.\tag4.33 
$$
\endproclaim 

Note in particular that
$$\align
\eta _{C,V'_0}&=\eta (C,0)+\operatorname{dim}V_0(C),\quad \text{ if
}\Pi _1=\Pi 
_\ge(C),\tag4.34\\
\eta _{C,V'_0}&=\eta (C,0),\quad \text{ if }\operatorname{dim}V'_0=\tfrac12 \operatorname{dim}V_0(C).\tag4.35
\endalign$$

Note also that since $\zeta
(P_T,0)=a''_0(I)-\operatorname{dim}V_0(P_T)$, we have in 
the situation of Theorem 4.9 that
$$
\zeta (P_T,0)=-\tfrac14\eta _{C,V'_0}-\operatorname{dim}V_0(P_T)+
\text{ local contributions.}\tag4.36
$$

Under special circumstances, we can show that $\tilde c''_0$ and
$\tilde d''_0$ are purely local:

\proclaim{Theorem 4.11} Let $\Pi _1$ be an orthogonal
pseudodifferential projection (so that $\Pi _2=\Pi _1^\perp$), and
assume that there exists a unitary 
morphism
$\sigma $ such that $$
\sigma ^2=-I,\quad \sigma P'=P'\sigma ,\quad \Pi _1^\perp=-\sigma \Pi
_1\sigma .\tag4.37
$$
Then$$\aligned
\Tr_{X'}(\tfrac12\Pi _2\tfrac{\partial _\lambda
^{m-1}}{(m-1)!}(P'-\lambda )^{-1})&= \tfrac14\Tr_{X'}(\tfrac{\partial _\lambda
^{m-1}}{(m-1)!}(P'-\lambda )^{-1}),\\
-\Tr_{X'}(\tfrac12\sigma \Pi _2\tfrac{\partial _\lambda
^{m-1}}{(m-1)!}(P'-\lambda )^{-\frac12})&= -\tfrac14\Tr_{X'}(\sigma
\tfrac{\partial _\lambda 
^{m-1}}{(m-1)!}(P'-\lambda )^{-\frac12}).
\endaligned\tag4.38$$
Thus in {\rm (4.8)}, 
$\tilde c''_0(I)$ ($=c''_0(I)$) is locally determined (from the
symbol of $P'$), and in {\rm
(4.14)} with $\psi =\sigma $,$$
\tilde d'_0(\sigma )=0, \quad \tilde d''_0(\sigma )\text{ is locally determined}\tag4.39
$$
(from the symbol of $P'$ and $\sigma $).

It follows that in the situation of Theorems {\rm 4.2} and 
{\rm 4.3} with $D_1=\sigma (\partial _{x_n}+B_1)$, $$\gathered
\tilde
a'_0(D_1)=a'_0(D_1)=0;\\
\tilde a''_0(I ),a''_0(I)\text{ and }\tilde a''_0(D_1)\text{ are
locally determined}\endgathered\tag4.40
$$
(depending only on finitely many homogeneous terms in the symbols of $P$ and $T$,
resp\. $P$, $T$ and $D_1$).

\endproclaim 
\demo{Proof}
Introduce the shorter notation $$
R_{m,1}=\tfrac{\partial _\lambda
^{m-1}}{(m-1)!}(P'-\lambda )^{-1},\quad R_{m,2}=\tfrac{\partial _\lambda
^{m-1}}{(m-1)!}(P'-\lambda )^{-\frac12};
$$ note that $\sigma R_{m,i}=R_{m,i}\sigma $ by (4.37). Then we have
for the traces on $X'$, using moreover that the trace is invariant
under circular permutation:
$$
\multline
\Tr(\Pi _2R_{m,1})=
-\Tr(\sigma \Pi _1\sigma R_{m,1})=
-\Tr( \Pi _1 R_{m,1}\sigma ^2)\\
=\Tr( (1-\Pi _2) R_{m,1} )=
\Tr( R_{m,1})-\Tr( \Pi _2 R_{m,1} ).\endmultline \tag4.41
$$
 It follows that$$
\Tr(\Pi _2R_{m,1})=\tfrac12 \Tr( R_{m,1}).\tag 4.42
$$ 
Similarly,
$$
\multline
\Tr(\sigma \Pi _2R_{m,2})=
-\Tr(\sigma ^2\Pi _1\sigma R_{m,2})=
\Tr( \Pi _1 \sigma R_{m,2})\\
=\Tr( (1-\Pi _2) \sigma R_{m,2} )=
\Tr(\sigma  R_{m,2})-\Tr(\sigma  \Pi _2 R_{m,2} ),\endmultline \tag4.43
$$implying
$$
\Tr(\sigma \Pi _2R_{m,2})=\tfrac12 \Tr(\sigma  R_{m,2}).\tag 4.44
$$ 
This shows (4.38).
It is classically known that  $\Tr(
R_{m,1})$ has an expansion in powers with
local coefficients; this shows the assertion on
$\tilde c''_0(I)=c''_0(I)$. For  $i=2$, $\Tr(
\sigma R_{m,2})$ has an expansion in powers of $(-\lambda )$ with only local
coefficients, since $\sigma R_{m,2}$
is strongly polyhomogeneous, cf\. \cite{GS95} or \cite{G02, Th\. 2.10}. This implies the
assertions on $\tilde d'_0(\sigma )$ and 
$\tilde d''_0(\sigma )$.

The last assertion now follows from (4.9) and (4.15). 
\qed\enddemo 

The proof also shows that $\tilde c'_0(I)=\tilde a'_0(I)=0$, but we
know that already from Theorem 4.5.

\subhead 5. Consequences for the APS problem \endsubhead

The preceding results have
interesting new consequences 
for the realizations of first-order operators in Example 1.2, which
we now consider in detail. Let $D$ satisfy (1.6) and let $\Pi $ be a
well-posed orthogonal $\psi $do projection for $D$. Then in view of
Green's formula:$$
(Du,v)_{X}-(u,D^*v)_X=-(\sigma \gamma _0u,\gamma _0v)_{X'},\text{ for
}u\in C^\infty (E_1),v\in C^\infty (E_2),\tag5.1
$$
the adjoint $(D_\Pi )^*$ is the realization of $D^*$ (of the form
$(-\partial _{x_n}+A+x_nA_{21}+A_{20})\sigma ^*$ on $X_c$) defined by the
boundary condition $\Pi ^\perp\sigma ^*\gamma _0v=0$ (associated with
the well-posed projection $\Pi '=\sigma \Pi ^\perp\sigma ^*$ for $D^*$).
It follows that $D^*D$ is of the form (1.1) with $P'=A^2$, and that
${D_\Pi }^*D_{\Pi }$ is the realization of $D^*D$ defined by the
boundary condition $$
\Pi \gamma _0u=0,\quad \Pi ^\perp(\gamma _1u+A_1(0)\gamma _0u)=0.\tag5.2
$$ 
Thus  ${D_\Pi }^*D_{\Pi }$ is of the type $P_T$ considered in
Sections 1--4, with $P=D^*D$, $P'=A^2$,
 $\Pi _1=\Pi $ and $B=A_1(0)$. Note that the symbol considered in Assumption
2.7 is here$$
((a^0)^2+\mu ^2)^{\frac12}-(I-\pi ^h)a^0(I-\pi ^h).\tag 5.3
$$

When the principal symbols of $\Pi $ and $A^2$ commute, Assumption
2.7 is essentially equivalent with well-posedness. More precisely, we
have:

\proclaim{Lemma 5.1} Let $\Pi $ be an orthogonal $\psi $do projection
in $L_2(E'_1)$ and let $P_T$ be the realization of $D^*D$ under the
boundary condition {\rm (5.2)}, and assume that the principal symbols
of $\Pi $ and $A^2$ commute.

$1^\circ$ When $\Pi $ is well-posed for $D$, then Assumption {\rm 2.7}
holds for $\{P,T\}$ with $\theta =\frac\pi 2$.

$2^\circ$ If Assumption {\rm 2.7} holds for $\{P,T\}$ with some
$\theta >0$ and 
$\pi ^0(x',\xi ')$ has rank $N/2$, then $\Pi $ is well-posed for $D$. 
\endproclaim 

\demo{Proof} $1^\circ$. Fix $x'$, $|\xi '|\ge 1$, and consider the
model realization $d^0_{\pi ^0}$ (defined for the ordinary
differential operator $d^0=\sigma (x')(\partial _{x_n}+a^0(x',\xi
'))$ in $L_2(\Bbb R_+,\Bbb C^N)$ by the boundary condition $\pi
^0(x',\xi ')u(0)=0$), and the model realization $p^0_{t^0}$ (defined
similarly from principal symbols). The well-posedness assures that
$d^0_{\pi ^0}$ is injective, hence $p^0_{t^0}=(d^0_{\pi
^0})^*d^0_{\pi ^0}$ is selfadjoint positive, as an unbounded operator
in  $L_2(\Bbb R_+,\Bbb C^N)$. It follows that $p^0_{t^0}-\lambda $ is
bijective from its domain to $L_2(\Bbb R_+,\Bbb C^N)$, for all
$\lambda \in\Bbb C\setminus \Bbb R_+$. Using that $\pi ^0$
commutes with $(a^0)^2$, we can carry out the calculations in
the proof of Lemma 2.6 for the model problem (without commutation error
terms), which 
allows us to conclude that the equation in
$\Bbb C^N$: $$
[((a^0)^2-\lambda )^{\frac12}-(I-\pi ^0)a^0(I-\pi ^0)]\varphi =\psi ,
\tag5.4$$
is uniquely solvable for $\psi \in R(I-\pi ^0)$. Moreover,
the calculations in
Remark 2.9 on the model level extend the solvability of (5.4) to all $\psi
\in\Bbb C^N$. 
The invertibility property extends readily to the strictly
homogeneous symbols for $\xi '\ne 0$, it is obvious for $\xi '=0$
with $\lambda \ne 0$.

$2^\circ$. Assumption 2.7 gives for $\mu =0$, $|\xi '|\ge 1$, that
$p^0_{t^0}=(d^0_{\pi ^0})^*d^0_{\pi ^0}$ is bijective. This implies
injectiveness of $d^0_{\pi ^0}$, i.e., injective ellipticity of
$\{d^0, \pi ^0\gamma _0\}$. Then well-posedness holds exactly when
$\pi ^0$ has rank $N/2$. (One may consult \cite{G99, p\. 55}.) 
\enddemo 

\example{Example 5.2} When $\Pi $ is taken as $\Pi _\ge+\Cal S$ with
$\Cal S$ of order $-1$ (cf\. Example 1.2),  $\pi ^0$ commutes with
$a^0$ itself, and we see directly that 
Assumption 2.7 holds simply
because \linebreak $-(I-\pi ^h)a^0(I-\pi
^h)\ge 0$. ---
For the projections $\Pi (\theta )=P(\theta )$ introduced
by Br\"uning and 
Lesch in \cite{BL99}, Assumption 2.7 is also directly verifiable, since the conditions of 
\cite{BL99} assure that $-(I-\pi ^h)a^0(I-\pi ^h)- c |a^0 |\ge 0$
for some $c>-1$. Here $\Pi (\theta )$ commutes with $A^2$. Again,
perturbations of order $-1$ are allowed.
\endexample

Thus the results of Section 2--4 apply to ${D_\Pi }^*D_{\Pi }$ with
$\Gamma _\theta =\Gamma $.
So there are expansions (2.43), (3.4) and (3.9) for$$
\gathered
\Tr(\varphi ({D_\Pi }^*D_\Pi -\lambda )^{-m}),\quad 
\Tr(\varphi e^{-t{D_\Pi }^*D_\Pi}),\\
\text{ and }\Gamma (s)\Tr(\varphi ({D_\Pi }^*D_\Pi)^{-s})=\Gamma (s)\zeta
(\varphi , {D_\Pi }^*D_\Pi,s),\endgathered
$$
and, with the choice $F=\varrho D$ ($\varrho $ a morphism from $E_2$
to $E_1$), there are expansions as in (2.35), (3.5) and (3.11)--(3.12) for$$
\gathered
\Tr(\varrho D ({D_\Pi }^*D_\Pi -\lambda )^{-m}),\quad 
\Tr(\varrho D e^{-t{D_\Pi }^*D_\Pi}),\quad
\Gamma (s)\Tr(\varrho D ({D_\Pi }^*D_\Pi)^{-s}),\\
\text{ and
}\Gamma (\tfrac {s+1}2)\Tr(\varrho D ({D_\Pi }^*D_\Pi)^{-\frac{s+1}2})
=\Gamma (\tfrac{s+1}2)\eta
(\varrho  , D_\Pi ,s).\endgathered
$$
Such expansions were shown in \cite{G99} by a different procedure where
$D$ was regarded as part of a first-order system of the double size.

We get new results by drawing some consequences for the
coefficients at $k=0$ from Section 4.

Before doing this, let us also briefly look at $D_\Pi {D_\Pi }^*$. It
is easily checked 
that $\sigma 
^*DD^*\sigma $ is of the form (1.1), and that $\sigma ^*D_\Pi {D_\Pi
}^*\sigma $ is the realization of it with boundary condition$$
\Pi ^\perp\gamma _0u=0,\quad \Pi (\gamma _1u-(A+A_{20}(0))\gamma _0u)=0.\tag5.5 
$$
In the consideration of trace formulas for $D_\Pi {D_\Pi }^*$, a
composition to the left with $\sigma $ and to the right with $\sigma
^*$ leaves the formulas corresponding to 
(2.43), (3.4) and (3.9) 
unchanged if $\varphi =I$.

Theorems 4.5 and 4.7 imply immediately:

\proclaim{Corollary 5.3} Let $P_T=D_\Pi ^*D_\Pi $, where $D$ is as in
Example {\rm 1.2}, $\Pi $ is well-posed for $D$, and the principal
symbols of $\Pi $ and $A^2$ commute.

\flushpar{\rm (i)} For the expansions {\rm (2.43), (3.4), (3.9),}
related to the zeta function,$$
\tilde a'_0(I)=a'_0(I)=0.\tag5.6
$$ 
Moreover,$$
\tilde a'_0(\varphi )=a'_0(\varphi
)=\tfrac14\operatorname{res}(\varphi \Pi ^\perp);\tag5.7
$$
it is zero in the following cases {\rm (a)} and {\rm (b)}:
\roster
\item"{\rm (a)}" $\varphi \Pi ^\perp$ is a projection,
\item"{\rm (b)}" $n$ is odd and $\Pi ^\perp=\Pi _>(C)+\Cal S$ for some first-order
selfadjoint elliptic {\bf differential} operator $C$ of order $1$, $\Cal S$ of
order $-n$.
\endroster

\flushpar{\rm (ii)} For the expansions {\rm (2.35), (3.5),
(3.11)--(3.12)}
with $F=\varrho D$,
related to the eta function,
$$
\tilde a'_0(\varrho D )=a'_0(\varrho D)=\alpha
\operatorname{res}(\varrho \sigma  \Pi ^\perp);\tag5.8
$$
it is zero if $\varrho \sigma \Pi ^\perp$ is a projection, or if 
{\rm (b)} holds.
\endproclaim   

Note that (5.6) means that the zeta function of $D_\Pi
^*D_\Pi $ is regular at zero. Since the hypotheses assuring this
(once $D$ is taken of the form (1.6)), are entirely concerned with
principal symbols, we have in particular: {\it The regularity of the
zeta function at $s=0$ is preserved under perturbations of $\Pi $ of
order $-1$}. When $[\pi ^0,(a^0)^2]=0$, this is a far better result
than that of \cite{G01$'$}, 
where it was shown for perturbations of order $-n$.

We also have from Theorem 4.9 and the following considerations:

\proclaim{Corollary 5.4} Assumptions as in Corollary {\rm 5.3}. Let
$\Pi $ be a spectral projection as in Definition {\rm 4.6}, with the
notation introduced there and in Definition {\rm 4.10}. Then in the
expansions {\rm (2.43), (3.4), (3.9),} $$\gather
 a''_0(I)=-\tfrac14 \eta _{C,V'_0}+\text{ local contributions,}
\tag 5.9\\
\zeta (D_\Pi ^*D_\Pi ,0)=-\tfrac14\eta
_{C,V'_0}-\operatorname{dim}V_0(D_\Pi )+\text{ local contributions.}
\tag5.10
\endgather$$
\endproclaim 

There is a similar result for $D_\Pi D_\Pi ^*$; here $\Pi $ is
replaced by $\Pi ^\perp=\Pi _>(-C)+\Pi _{V''_0}$ in view of the
remarks above on $D_\Pi D_\Pi ^*$. So in this case, Theorem 4.9 gives:$$
\aligned
 a''_0(I)(D_\Pi D_\Pi ^*)&=-\tfrac14 (-\eta
(C,0)+\operatorname{dim}V''_0-\operatorname{dim}V'_0)
+\text{ local contributions}\\
&=\tfrac14 \eta _{C,V'_0}+\text{ local contributions,} 
\endaligned\tag 5.11
$$
 
Observe moreover that since 
$$
\operatorname{index}D_{\Pi }=\Tr e^{-tD_\Pi ^*D_\Pi }-\Tr e^{-tD_\Pi
D_\Pi ^*}=
a''_0(I)(D_\Pi ^*D_\Pi )-a''_0(I)(D_\Pi D_\Pi^*),\tag5.12
$$ 
we find: 

\proclaim{Corollary 5.5} In the situation of Corollary {\rm 5.4},$$
\operatorname{index}D_{\Pi }=-\tfrac12\eta _{C,V'_0}+ \text{ local
contributions.} \tag 5.13
$$
\endproclaim  

For the case where $\Pi =\Pi _\ge(A)$ (cf\. (4.34)), (5.13) is known
from \cite{APS75}, 
and (5.9) is known from \cite{G92}; for $\Pi =\Pi _>(A)+\Pi
_{V'_0}$ in the product case, cf\. \cite{GS96, Cor\. 3.7}.
We believe that it is
an interesting new result that for rather general projections,
the non-locality 
depends only on the projection, not the interior operator, in this sense.

Now we turn to cases with selfadjointness properties. 
We are here both interested in truly selfadjoint product cases and
in nonproduct cases where $D$ is principally
selfadjoint at $X'$.
Along with $D$ we consider the operator of product type $D_0$ defined
by
$$
D_0=\sigma (\partial _{x_n}+A)\text{ on }X_c, \text{ so that
}D=D_0+\sigma (x_nA_{11}+A_{10}).\tag5.14
$$ 
In addition to
the requirements that $\sigma $ be unitary and $A$ be selfadjoint, we
now assume that $E_1=E_2$ and that $D_0$ is formally selfadjoint on
$X_c$ when this is provided with the ``product'' volume element
$v(x',0)dx'dx_n$; 
this means 
that$$
\sigma ^2=-I,\quad \sigma A=-A\sigma .\tag 5.15
$$
$D_0$ can always be extended to an elliptic operator on $X$ (e.g\. by
use of $D$); let us denote the extension $D_0$ also. If the extension
is selfadjoint, we call this a {\it selfadjoint product case}.

When $\Pi $ is an orthogonal projection in $L_2(E'_1)=L_2(E'_2)$, it
is well-posed for $D$ if and only if it is so for $D_0$.
For $D_0$ in selfadjoint product cases,
some choices of $\Pi $ will lead to selfadjoint realizations
$D_{0,\Pi }$, namely 
(in view of (5.1)) those for which
$$
\Pi =-\sigma \Pi ^\perp \sigma .\tag 5.16
$$

The properties (5.15) and (5.16) imply (4.37) with $P'=A^2$, $\Pi
_1=\Pi $, so we can apply Theorem
4.11 to  $D_\Pi ^*D_{\Pi }$ (and $D_{0,\Pi }^2$).

As pointed out in the appendix A.1 of Douglas and Wojciechowski
\cite{DW91}, it follows from Ch\. 17 (by Palais and Seeley) of Palais \cite{P65}
that when (5.15) holds and $n$ is odd, there exists a subspace $L$
of $V_0(A)$ such that $\sigma L\perp L$ and $V_0(A)= L\oplus\sigma
L$. M\"uller showed in \cite{M94} (cf\. (1.5)ff\. and Prop\.
4.26 there) that such $L$ can be found in any dimension. Denoting the
orthogonal projection onto $L$ by $\Pi _L$, we have that  
$$
\Pi _+=\Pi _>(A)+\Pi _L\tag5.17
$$
satisfies (5.16). The projections $\Pi (\theta )$
introduced by Br\"uning and Lesch \cite{BL99} likewise satisfy (5.16). 

We here conclude from Theorem 4.11:

\proclaim{Corollary 5.6} In addition to the assumptions of Corollary
{\rm 5.3}, assume that $E_1=E_2$ and that {\rm (5.15), (5.16)} hold.
Then in 
{\rm (2.43), (3.4), (3.9)} 
for $P_T=D_\Pi ^*D_{\Pi }$, $\tilde a''_0(I)$
($=a''_0(I)$) is locally determined. Equivalently,
the value of $\zeta (D_\Pi ^*D_\Pi ,s)$ at
zero satisfies:$$
\zeta (D_\Pi ^*D_\Pi ,0)=-\operatorname{dim}V_0(D_\Pi )+\text{ local contributions}.\tag5.18
$$
\endproclaim 

We use this to show for the zeta function:

\proclaim{Theorem 5.7} In addition to the hypotheses of Corollary
{\rm 5.6}, assume that$$
\Pi =\overline\Pi +\Cal S,\tag5.19
$$
where $\overline\Pi $ is a fixed well-posed projection satisfying {\rm (5.16)} and $\Cal S$ is of order $\le
-n$. ($\overline \Pi $ can in particular be taken as $\Pi _+$ in {\rm
(5.17)} or $\Pi (\theta )$ from {\rm \cite{BL99}}.)

Then the $\tilde a''_0$-terms (and $a''_0$-terms) in {\rm (2.43),
(3.4), (3.9)} for $D_{\overline\Pi }$ and 
$D_\Pi $ are the same,$$
\tilde a''_0(I)(D_\Pi ^*D_\Pi )=\tilde a''_0(I)
(D_{\overline\Pi }^*D_{\overline\Pi }),\tag5.20
$$
so
$$
\zeta (D_\Pi ^*D_\Pi  ,0)+\operatorname{dim}V_0(D_\Pi )=\zeta 
(D_{\overline\Pi }^*D_{\overline\Pi },0)+\operatorname{dim}V_0(D_{\overline\Pi } ),\tag5.21
$$
and in particular,$$
\zeta (D_\Pi ^*D_\Pi ,0)=\zeta (D_{\overline\Pi }^*D_{\overline\Pi },0) \quad (\operatorname{mod} \Bbb Z).\tag5.22
$$
\endproclaim 

\demo{Proof}
We shall combine the fact that $\tilde a''_0(I)$ is locally
determined with order considerations.  Let $$ 
R_{\overline T}^m=(D_{\overline\Pi }^*D_{\overline\Pi }-\lambda )^{-m},\quad R_{T}^m=(D_\Pi
^*D_\Pi 
-\lambda )^{-m}.\tag5.23
$$
Note that they have the same pseudodifferential part $(D^*D-\lambda
)_+^{-m}$, so their difference $R_{\overline T}^m-R_{T}^m$ is a singular
Green operator. It is shown in \cite{G01$'$, proof of Th\. 1} that when $\Cal
S$ is of order $ -n$, the $\psi $do $\tr_n(R_{\overline T}^m-R_{T}^m)$
on $X'$ has symbol in $S^{-m-n,-m,0}\cap S^{-n,-2m,0}$. The total
order is  $-n-2m$, so the highest degree of the homogeneous terms in
the symbol is $-n-2m$. As noted in Remark 4.4, the local contribution to the
terms with index $k=0$ in the trace expansion of this difference
comes from homogeneous terms of degree $1-n-2m$ (recall that $\operatorname{dim}X'=n-1$), so
since the terms consist purely of local contributions, they must
vanish. This shows 
(5.20), and the other consequences are immediate.
\qed

\enddemo 

In a language frequently used in this connection, the statement means that 
$\zeta (D_\Pi ^*D_\Pi  ,0)+\operatorname{dim}V_0(D_\Pi )$ {\it is constant on the 
Grassmannian of $\psi $do projections satisfying {\rm (5.16)} and differing
from $\overline\Pi $ by a term of order}
$\le -\operatorname{dim}X$.

When $D$ equals $D_0$ in a selfadjoint product case,
Wojciechowski shows a 
result like this in \cite{W99, Sect\. 3} for 
perturbations of $\Pi _+$ of
order $-\infty $, assuming that $D_{0,\Pi _+}$ and
$D_{0,\Pi }$ are  
invertible. The non-invertible case is
treated by Y\. Lee in the appendix of
\cite{PW02}; he shows moreover that $\zeta (D^2_{0,\Pi
_+},0)+\operatorname{dim}V_0(D_{0,\Pi _+} )=0$, so we conclude that 
$\zeta (D_{0,\Pi }^2,0)+\operatorname{dim}V_0(D_{0,\Pi } )=0$ when
$\Pi =\Pi _++\Cal S$. 
\medskip

We can also discuss the eta function $\eta (D_\Pi ,s)=\Tr(D(D^*_\Pi
D_\Pi )^{-\frac{s+1}2})$, extended meromorphically as in (3.12).
First we conclude from Theorem 4.11:

\proclaim{Corollary 5.8} Assumptions of Corollary
{\rm 5.6}.
In
{\rm (2.35), (3.5), (3.11)--(3.12)} for $P_T=D_\Pi ^*D_\Pi $, $D_1=D$, 
one has that $\tilde a'_0(D)=a'_0(D)=0$, and 
$\tilde a''_0(D)$
($=a''_0(D)$) is locally determined. 

In other words, the double pole
of $\eta (D_\Pi ,s)$ at $0$ vanishes and the residue at $0$ is
locally determined.

\endproclaim 

(It may be observed that when (5.15) holds, the entries in the second
line of (4.38) vanish identically, since $\Tr_{X'}(\sigma \partial _\lambda
^{m-1}(P'-\lambda )^{-\frac12})=0$.)

It is remarkable here that the hypotheses, besides  
(5.15)--(5.16), only contain
requirements on principal symbols (the well-posedness of $\Pi $ for
$D$ and the
commutativity of the principal symbols of $\Pi $ and $A^2$). So the
result implies in particular that {\it the vanishing of the double pole of
the eta function is invariant under perturbations of $\Pi $ of order
$-1$} (respecting (5.16)). Earlier results have dealt with
perturbations of order $-\infty $ \cite{W99} or order $-n$ \cite{G01$'$}.

Now consider the simple pole of $\eta (D_\Pi ,s)$ at 0. Here we can
generalize the result of Wojciechowski 
\cite{W99} on the 
regularity of the eta function after a
perturbation of order $-\infty $, to perturbations of order $-n$ of
general $\overline\Pi $:

\proclaim{Theorem 5.9} Assumptions of Theorem {\rm 5.7}.

In {\rm (2.35), (3.5), (3.11)--(3.12)} with
$D_1=D$, 
the $\tilde a''_0$-terms (and $a''_0$-terms) for $D_{\overline\Pi }^*D_{\overline\Pi }$ and
$D_\Pi ^*D_\Pi $ are the same:$$
\tilde a''_0(D)(D_\Pi ^*D_\Pi )=\tilde a''_0(D)(D_{\overline\Pi }^*D_{\overline\Pi });\tag5.24
$$
in other words, $\operatorname{Res}_{s=0}\eta (D_\Pi
,s)=\operatorname{Res}_{s=0}\eta (D_{\overline\Pi },s)$.

In particular, if $\tilde a''_0(D)(D_{\overline\Pi }^*D_{\overline\Pi
})=0$ (this 
holds for $\Pi _+$ and for certain $\Pi (\theta )$ if $D$ equals
$D_0$ in a selfadjoint product case), then $\tilde a''_0(D)(D_\Pi
^*D_\Pi )=0$, i.e., 
the eta function 
$\eta (D_{\Pi },s)$ is regular at $0$.

\endproclaim 

\demo{Proof}
As in the proof of Theorem 5.7,
we  combine the fact that $\tilde a''_0(D)$ is locally
determined with order considerations.  Consider $DR_{\overline T}^m$ and 
$DR_{T}^m$, cf\. (5.23).
Since they have the same pseudodifferential part $D(D^*D-\lambda
)_+^{-m}$, their difference $DR_{\overline T}^m-DR_{T}^m$ is a singular
Green operator. It is shown in \cite{G01$'$, proof of Th\. 1} that when $\Cal
S$ is of order $ -n$, the $\psi $do $\tr_n(DR_{\overline T}^m-DR_{T}^m)$
on $X'$ has symbol in $S^{-m-n,1-m,0}\cap S^{-n,1-2m,0}$. The total
order is  $1-n-2m$, so the highest degree of the homogeneous terms in
the symbol is $1-n-2m$. Now the local contribution to the
terms with index $k=0$ in the trace expansion of this difference
comes from homogeneous terms of degree $2-n-2m$ (cf\. Remark 4.4), so
since the terms contain only local contributions, they must vanish. This shows
(5.24).

In particular, $\tilde a''_0(D)(D_{\Pi }^*D_\Pi )$ vanishes if $\tilde
a''_0(D)(D_{\overline\Pi }^*D_{\overline\Pi })$ does so; then the eta
function for 
$D_{\Pi }$ is
regular at 0.
The eta regularity for the case $\overline\Pi=\Pi _+$, $D$ equal to $D_0$ and
selfadjoint on $X$ with product volume element on $X_c$, was shown in
\cite{DW91} under 
the assumptions $n$ odd and $D$ compatible; this was extended to
general $n$ and not necessarily compatible $D$ in M\"uller \cite{M94}. It
was shown for 
certain $\Pi (\theta )$ in \cite{BL99, Th\. 3.12}. 
\qed

\enddemo

In a frequently used terminology, the theorem shows that the residue
of the eta function 
is constant on the 
Grassmannian of $\psi $do projections satisfying (5.16) and differing
from $\overline\Pi $ by a term of order
$\le -\operatorname{dim}X$. We do not expect that the
order can be lifted further in general.

The result on the regularity of the eta function at $s=0$ for
$(-n)$-order perturbations of the product case with $\overline \Pi
=\Pi _+$ has 
been obtained independently by Yue Lei \cite{L02} at the same time as our
result, by another
analysis based on heat operator formulas. 

The above results on the vanishing of the eta residue are concerned with
situations where $D$ equals $D_0$ in a selfadjoint
product case. However, the fact from Corollary 5.8
that $\operatorname{Res}_{s=0}\eta (D_\Pi ,s)$ is locally determined
also in suitable non-product cases, should facilitate the
calculation of the residue then. For example, if $D=D_0+x_n^{n+1}A_{31}$ on
$X_c$ with a first-order tangential differential operator $A_{31}$,
and the volume element satisfies $\partial _{x_n}^jv(x',0)=0$ for
$1\le j\le n$, then by
\cite{G02, proof of Th\. 3.11}, the local terms with index $k\le 0$ in the
difference between the resolvent powers are determined entirely
from the interior operators $D$ and $D_0$. Then when $n$ is even, the
contributions to $k=0$ vanish simply because of odd parity in $\xi $;
this gives examples where the eta function is regular at $0$ in a
non-product situation.

\enddocument